\documentclass[3p,times]{elsarticle}

\usepackage{graphicx}
\usepackage{psfrag}
\usepackage[dvipsnames]{xcolor}

\usepackage{amsmath}
\usepackage{amsthm}
\usepackage{amsfonts}
\usepackage{dsfont}
\usepackage{rotating}
\usepackage{multirow}
\usepackage{epstopdf}
\usepackage{newtxmath}
\usepackage{mathtools}
\usepackage{bm}
\usepackage{bbm}
\usepackage{comment}
\usepackage{booktabs}
\usepackage{accents}
\usepackage{empheq}

\usepackage{tikz}
\usepackage{pgfplots}
\usepackage{pgffor}
\usepackage{booktabs}
\usepackage{rotating}
\usepackage{lscape}
\usepackage{graphics}

\newtheorem*{remark}{Remark}
\newtheorem{definition}{Definition}

\DeclareMathAlphabet\mathbfcal{OMS}{cmsy}{b}{n}  

\newcommand{\CFL}{\textnormal{CFL}}
\newcommand{\cc}{\mathbf{c}}
\newcommand{\cell}[1]{\omega_#1}

\newcommand{\dt}{{\Delta t}}

\newcommand{\F}{\mathbf{F}}

\newcommand{\f}{\mathbf{f}}

\newcommand{\HH}{\mathcal{H}}

\newcommand{\I}{\mathbf{I}}

\newcommand{\nn}{\mathbf{n}}
\newcommand{\Omegabg}{\Omega_{\text{bg}}}
\newcommand{\Omegafg}{\Omega_{\text{fg}}}

\newcommand{\QQ}{\mathbf{Q}}

\newcommand{\R}{\mathds{R}}

\newcommand{\Ste}[1]{\mathcal{S}_#1}
\newcommand{\Tbg}{\mathcal{T}_{\text{bg}}}
\newcommand{\Tfg}{\mathcal{T}_{\text{fg}}}
\newcommand{\TOmega}{\mathcal{T}_{\Omega}}

\newcommand{\uu}{\mathbf{u}}
\newcommand{\w}{\mathbf{w}}
\newcommand{\xx}{\mathbf{x}}

\newcommand{\zz}{\mathbf{z}}

\newcommand{\IM}[1]{#1_{\text{I}}}
\newcommand{\EX}[1]{#1_{\text{E}}}
\newcommand{\old}[1]{#1_{\text{old}}}

\newcommand{\Pol}[1]{\mathcal{R}^P_#1}
\newcommand{\Polb}[1]{\mathcal{R}^Q_{#1}}

\newcommand{\diff}{\text{d}}
\newcommand{\Rey}{\text{Re}}
\newcommand{\ustar}{\uu^*}
\newcommand{\Tstress}{\mathbf{T}}
\newcommand{\St}{\text{St}}

\begin{document}
	
	\begin{frontmatter}
		
		\journal{Journal of Computational Physics}
		
		

		\title{Arbitrary-Lagrangian-Eulerian finite volume IMEX schemes for the incompressible Navier-Stokes equations on evolving Chimera meshes}
		
        \author[dmi]{Michele Giuliano Carlino}
        \ead{michelegiuliano.carlino@unife.it}
        
        \author[dmi]{Walter Boscheri$^*$}
        \ead{walter.boscheri@unife.it}
        \cortext[cor1]{Corresponding authors}
        
		\address[dmi]{Department of Mathematics and Computer Science, University of Ferrara, Via Machiavelli 30, 44121 Ferrara, Italy}
%
\begin{abstract}
In this article we design a finite volume semi-implicit IMEX scheme for the incompressible Navier-Stokes equations on evolving Chimera meshes. We employ a time discretization technique that separates explicit and implicit terms, accommodating the multi-scale nature of the governing equations, which encompass both slow and fast scales. The finite volume approach for both explicit and implicit terms allows to encode into the nonlinear flux the velocity of displacement of the Chimera mesh via integration on moving cells. The numerical solution is then projected onto the physically meaningful solution manifold of non-solenoidal fields that stems from the energy equation. To attain second-order time accuracy, we employ semi-implicit IMEX Runge-Kutta schemes. These novel schemes are combined with a fractional-step method, thus the governing equations are eventually solved using a projection method to satisfy the divergence-free constraint of the velocity field. The implicit discretization of the viscous terms allows the CFL-type stability condition for the maximum admissible time step to be only defined by the relative fluid velocity referred to the movement of the frame and not depending also on the viscous eigenvalues. Communication between different grid blocks is enabled through compact exchange of information from the fringe cells of one mesh block to the field cells of the other block. Taking advantage of the continuity of the solution and the definition of a minimal compact stencil, the numerical solution of any system of differential equations is characterized by continuous data extrapolation. In this way, the continuity of the solution is recovered in one-shot during the solution of the arising algebraic systems by not involving neither direct discretization of the differential operators on fringe cells nor an iterative Schwartz-type method. Free-stream preservation property, i.e. compliance with the Geometric Conservation Law (GCL), is respected. The accuracy and capabilities of the new numerical schemes is proved through an extensive range of test cases, demonstrating ability to solve relevant benchmarks in the field of incompressible fluids.
\end{abstract}
%
\end{frontmatter}



\section{Introduction}
Incompressible flows are mathematically described by nonlinear systems of hyperbolic conservation laws, encompassing a wide array of physical phenomena such as environmental, geophysical, and meteorological flows, as well as the dynamics of mechanical processes like turbo-machinery or turbines in the energy engineering. Within this work, we focus on the incompressible Navier-Stokes equations, utilized to model atmospheric flows, flows around wings, and even pressurized pipe flows. Given the significance of these physical applications, there is considerable interest in developing numerical solutions for this physical model, which, however, poses certain complexities. In particular, the energy equation in the incompressible Navier-Stokes model reduces to a constraint on the velocity field, imposing a divergence-free condition. The adoption of unstructured grids proves highly advantages for real-world applications, accommodating ocean bathymetry, wing and turbine shapes, and river morphology, all of which require high-accuracy approximations of the geometry. Additionally, attaining high-order accuracy in time represents a crucial objective to obtain precise results for unsteady problems. In general, the simulation of complex flows with unsteady geometries (e.g., fluid-structure interactions, freely moving objects or moving boundaries induced by the moving flow itself) needs the modeling of \textit{ad hoc} schemes. The state of the art for this kind of methods is defined by three main approaches: i) Arbitrary-Lagrangian-Eulerian (ALE) methods, ii) interface approaches and iii) the employment of Chimera meshes.

The ALE methods \cite{hirt1997arbitrary, duarte2004arbitrary} are known for their accuracy and ability to handle complex grid displacement and mesh adaptation through a reformulation of the governing equations. However, when the grid deformation becomes excessively stretched or twisted, these schemes can lead to computationally expensive remeshing requirements. Strong differential rotations can also become problematic with this approach. Consequently, additional numerical errors may arise due to the interpolation of data from the old grid to the new mesh, necessitating careful management. To overcome the issues induced by the interpolation, direct ALE methods with topology changes have been recently forwarded in \cite{GaburroAREPO,Re2022}, but they still need to reinitialize the entire topology at each time step of the simulation.
In contrast, interface methods, like ghost-boundary methods \cite{gibou2002second,COCO2020109623}, immersed boundary methods \cite{mittal2005immersed}, and penalization methods \cite{angot1999penalization}, discretize the physical domain using a simple mesh, typically structured and Cartesian, which remains fixed throughout the simulation \cite{glowinski1994fictitious, peskin2002immersed}. This approach may not always perfectly fit the moving boundary, and special attention is required to achieve sufficient accuracy at the physical boundaries. Nevertheless, due to the simplicity and uniform aspect ratio of the mesh, the presence of thin boundary layers could significantly impact the computational advantages of these methods.
To address this concern, hybrid techniques have been developed, which combine immersed boundary methods with anisotropic mesh adaptations \cite{abgrall2014immersed}. This combination helps to overcome the challenges posed by thin boundary layers and it ensures a more efficient and accurate simulation.

Our investigations primarily focus on Chimera grids \cite{volkov1970method, benek19853, meakin1999composite, petersson1999hole}, which are composed of multiple overlapping mesh blocks that together form an overset grid \cite{starius1977constructing, starius1980composite, starius1977composite}. The Chimera meshing technique involves a discretization of the possibly evolving computational domain through a grid embedding approach. Initially, a major grid (background grid) is constructed, which is non-conformal with respect to the complex shape of the domain. Then, minor grids (foreground blocks) are created to describe the specific shape of regions containing obstacles. The minor blocks overlap with the major grid, establishing an overlapping region among all the blocks \cite{meakin1999composite}. This mesh generation strategy significantly simplifies mesh adaptation tasks, especially for scenarios involving boundary layers, changing geometries in unsteady problems, and unsteady multiply connected domains \cite{banks2016added, banks2013stable, schwendeman2010study, banks2007high, chesshire1990composite, henshaw2012cgins}.
Typically, numerical solutions on Chimera grids involve exchanging data through fringe cells located at the overlapping zone. One approach involves donor cells from a block in close proximity to the overlapping zone providing information to receptor cells of another block through polynomial interpolation, as proposed in \cite{wu2019numerical, guerrero2006overset, wang1995fully, zheng2003novel}. Another strategy, as presented in \cite{henshaw2005multigrid}, is concerned with the automatic generation of a coarse grid, and interpolation information is connected at the overlapping zone using a multigrid approach. Domain Decomposition (DD) methods, such as Schwartz, Dirichlet/Neumann, or Dirichlet/Robin methods, also enable communication between the different blocks. Here, each mesh block is considered as an independent domain decomposition, and the overlapping zones serve as interfaces for coupling the blocks. These approaches employ iterative discrete methods for two-way communication, and further details can be found in \cite{houzeaux2017domain}. Recently in \cite{bergmann2022second}, all involved local operators (such as Laplacian and normal gradient) are discretized also at the boundary of the communication interface by involving proper functional minimization in the sense of finite differences \cite{raeli2018finite} as there was only one block mesh.
There are also different approaches that connect the background and foreground meshes, such as the DRAGON grids \cite{kao1995advance}. DRAGON grids replace the overlapping zone with an unstructured grid during a subsequent stage while preserving the body-fitting advantages of the Chimera meshes. Essentially, a DRAGON grid generates a unified block mesh from a Chimera configuration. However, the computational costs associated with generating a DRAGON grid for an evolving domain can be significant, as a new DRAGON mesh must be built at every time instance.

From a numerical perspective, the multi-scale nature of the incompressible Navier-Stokes system imposes stringent limitations on the maximum allowable time step to ensure that the numerical scheme effectively captures fast waves. When using an explicit time discretization, this stability condition may yield exceedingly small time steps, rendering the method impractical for real-world applications. Additionally, the scheme introduces a significant amount of numerical dissipation, undermining the accuracy of the solution.
To overcome these challenges, flux splitting techniques constitute an effective strategy to separate the slow and fast scales of the physical problem. Explicit time discretization is used for the slow scales, typically associated with nonlinear convective terms related to the fluid velocity. Conversely, implicit time stepping is retained for the fast scales, which involve acoustic waves like the sound speed or the celerity. Furthermore, the incompressible Navier-Stokes viscous terms dictate a parabolic time step restriction, which also results in very small time steps for viscous-dominated flows. Therefore, an implicit treatment of these terms might be more convenient. This class of numerical methods falls under the category of implicit-explicit (IMEX) \cite{ascher1997implicit, boscarino2017asymptotic, pareschi2005implicit, boscarino2016high} or semi-implicit time schemes \cite{casulli1990semi, casulli1992semi, dumbser2013staggered}. In these methods, an algebraic system must be solved for the unknown physical quantity, which is discretized implicitly (e.g., pressure or velocity field in fluid dynamics). Typically, this system involves an elliptic equation that needs to be solved across the entire computational domain.
IMEX schemes have proven effectiveness in various contexts, including low Mach compressible flows \cite{park2005multiple, boscheri2021high, boscheri2022high}, magnetized plasma flows \cite{fambri2021novel}, pipe flows \cite{ioriatti2018semi, ioriatti2019posteriori}, and applications involving free-surface environmental flows \cite{tumolo2013semi, casulli1999semi} and atmospheric phenomena \cite{tumolo2015semi, orlando2023imex}. When using Cartesian meshes for the spatial discretization \cite{boscheri2021high, fambri2017semi}, finite difference schemes are often employed for the implicit terms, and higher-order extensions can be achieved straightforwardly by enlarging the stencil of the finite difference operators. Finite differences can also be utilized on orthogonal unstructured meshes like Voronoi tessellations \cite{boscheri2013semi, boscheri2019high, boscheri2020space} up to second order of accuracy, making the solution of the implicit algebraic system relatively simple.
However, when dealing with general unstructured meshes, solving the algebraic system becomes more challenging due to the complicated spatial discretization. Discontinuous Galerkin methods offer an elegant solution as they provide both high accuracy and compactness of the stencil \cite{tavelli2014high, tavelli2016staggered, orlando2022efficient}. On the other hand, finite volume methods can ensure conservation and excellent shock capturing properties. Therefore, they are suitable for discretizing explicit slow scale terms, like the nonlinear convective terms.
In order to merge the benefits of both methods, in the literature a new class of hybrid finite volume/finite element methods has emerged \cite{busto2018projection, bermudez2020staggered, busto2022staggered}. These hybrid methods combine the robustness of explicit finite volume solvers with the flexibility of the finite element method. However, these numerical methods were initially limited to simplex meshes in 2D/3D with triangles/tetrahedra. A recent effort in \cite{boscheri2023all} explored the usage of a hybrid scheme on general polygonal meshes, where the implicit pressure system is solved using a discontinuous Galerkin approach that acts on a sub-triangulation of the polygonal tessellation. Furthermore, in \cite{boscheri2023new} the coupling of different techniques involves the emergent Virtual Element Methods (VEM) \cite{beirao2013basic}. We address the reader to \cite{bergmann2022ader} for a hybrid finite volume/finite element methods for the incompressible Navier-Stokes equation on Chimera meshes. Explicit finite volume solvers on overset grids can be found in \cite{Ramirez2018}.

In this article, we introduce a finite volume (FV) IMEX scheme for the incompressible Navier-Stokes equations on evolving overset meshes. In particular, due to the movement of the mesh, the arising system of Partial Differential Equations (PDE) is no longer autonomous. Consequently, the discretization of the space-time overset configuration needs to properly address to the implicit as well as explicit physical variables during the integration of the governing equations. The FV solver is employed for treating the explicit terms. On the other hand, the integrated system along a generic control volume allows us to encode the displacement of the grid in the physical variables. This permits to solve the velocity and pressure field and, at the same time, to evolve the overset configuration as well as the computational domain. Finally, a second-order convergence is ensured by application of a semi-implicit IMEX Runge-Kutta time stepping method. Concerning the Chimera mesh, we provide an approach allowing the different blocks to communicate through a compact transmission \cite{bergmann2022second}. In particular, for a generic discretization of a differential problem, on fringe cells we force the solution to be directly extrapolated from the other partition without discretizing the differential operators and without employing any iterative process for collating the solution from different mesh blocks. Differently from \cite{bergmann2022ader}, in this work we use a method of lines in time, which eventually ensure free-stream preservation properties up to machine accuracy. Furthermore, the viscous sub-system is discretized implicitly by means of a compact finite volume scheme combined with a pressure-correction formulation. These methods first compute an intermediate velocity field that may not be solenoidal and then apply a correction to project the velocity onto the divergence-free manifold. We also refer to the recent work forwarded in \cite{MENG2020113040} where a fractional-step IMEX scheme is devised for incompressible flow with finite difference operators. In our work, we employ finite volume methods and different discretizations for the viscous terms as well as the class of semi-implicit IMEX schemes \cite{boscarino2016high} and no predictor-corrector scheme is adopted as done in \cite{MENG2020113040}.    

The paper is organized as follows. In Section \ref{sec.overset_grid}, the overset is introduced and the definition of minimal compact stencil is given. The mathematical model is presented in Section \ref{sec.math}. The numerical method is detailed in Section \ref{sec.numerical_method}. In particular, the discretization of all implicit and explicit terms is presented as well as the one-shot extrapolation allowing the different blocks to continuously exchange information in a compact transmission. The numerical results are shown in Section \ref{sec.numerical_results} for testing the convergence, robustness and accuracy of the novel numerical technique. Concluding remarks and possible future extensions ends the article in Section \ref{sec.concl}.

\section{Overset grid} \label{sec.overset_grid}
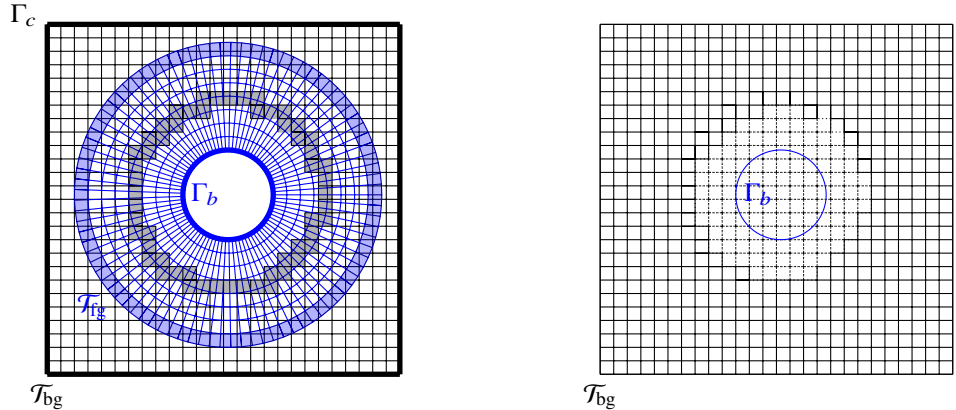
\begin{figure}
	\begin{center}
		\begin{tabular}{cc}
			\begin{tikzpicture}
				\begin{axis}[
					hide axis,
					axis equal,
					enlargelimits={abs=15pt},
					xlabel={x},
					ylabel={y},
					]
					
					\addplot [mark=none, black] table {overset_active_bg.dat};
					
					\addplot [mark=none, blue]  table {overset_active_fg.dat};
					
					
					\addplot [mark=none, fill=black, opacity=0.3]  table {fringe_bg.dat};
					
					\addplot [mark=none, fill=blue, opacity=0.3]  table {fringe_fg.dat};
					
					\addplot [mark=none, black, line width=2pt]  table {Gamma_c.dat};
					
					\addplot [mark=none, blue, line width=2pt]  table {Gamma_b.dat};
					
					\coordinate (Gammac) at (axis cs:-8,8);
				    \draw (Gammac) node[anchor=east]{$\Gamma_c$};

					\coordinate (Gammab) at (axis cs:-1,0);
					\draw (Gammab) node[blue]{$\Gamma_b$};

					\coordinate (Tbg) at (axis cs:-8,-8);
					\draw (Tbg) node[anchor=north]{$\Tbg$};

					\coordinate (Tfg) at (axis cs:-6,-4);
					\draw (Tfg) node[blue,anchor=north]{$\Tfg$};
					
				\end{axis}
				
			\end{tikzpicture} & 
		\begin{tikzpicture}
			\begin{axis}[
				hide axis,
				axis equal,
				enlargelimits={abs=15pt},
				xlabel={x},
				ylabel={y},
				]
				
				\addplot [mark=none, black] table {overset_active_bg.dat};
				
				
				\addplot [mark=none, black, densely dotted]  table {overset_hole.dat};
				
				
				
				
				\addplot [mark=none, black]  table {Gamma_c.dat};
				
				\addplot [mark=none, blue]  table {Gamma_b.dat};
				
				\coordinate (Gammab) at (axis cs:-1,0);
				\draw (Gammab) node[blue]{$\Gamma_b$};
				
				\coordinate (Tbg) at (axis cs:-8,-8);
				\draw (Tbg) node[anchor=north]{$\Tbg$};
				
			\end{axis}
			
		\end{tikzpicture}
		\end{tabular}
	\end{center}
	\caption{Sketch of an overset configuration (Chimera grid) for a domain $\Omega$ defined by a square $\Gamma_c$ without a circle  in the middle bounded by $\Gamma_b$. On the left, we plot all active cells of the background $\Tbg$ and foreground $\Tfg$ partitions. In particular, highlighted cells in black and blue denote fringe cells in background and foreground meshes, respectively. On the right, inactive cells in the hole are marked by dotted lines.}
	\label{fig.overset}
\end{figure}

The overset grid, also known as Chimera mesh, is a patch of different mesh blocks discretizing the computational domain $\Omega \subset \R^2$ (see Figure \ref{fig.overset}). Blocks overlap each other in a sub-region of the of the domain called \textit{overlapping zone}. In addition, the topology of any block can be different. In this work, all blocks are defined by quadrilateral cells. As explained in \cite{meakin1999composite}, firstly a Cartesian background mesh is built, and successively one or more foreground meshes are introduced inside the computational domain. The construction of the foreground blocks induces the definition of the overlapping regions, where the information deposited on the different meshes exchanges, and a coherent number of \textit{hole zones} in the background as depicted in Figure \ref{fig.overset} (right). In general, a total number $n_\ell$ of cell layers is imposed in both background and foreground for the overlapping zone. Consequently, the hole is defined by background cells completely covered by the foreground partition but not belonging to the appointed layers of superposition. 

This way of space discretization is useful when the boundary $\partial \Omega$ of the domain is composed of one or more internal boundaries $\Gamma_b$ and an external boundary $\Gamma_c$. For instance, in fluid dynamics this is the case of solid bodies $\Omega_b$, whose boundary $\partial \Omega_b$ coincides with $\Gamma_b$, immersed in a channel $\Omega_c$ of boundary $\Gamma_c$ and filled with a fluid. In this context, a Cartesian background grid is built regardless of the internal regions $\Omega_b$. Successively, the foreground meshes are employed in order to discretize the internal boundaries $\Gamma_b$. If regions $\Omega_b$ within the domain evolve over time, the foreground mesh will also evolve consistently with their displacement. This allows to avoid the remeshing of the different spatial configurations in time. Moreover, if there is a need to more accurately analyze the solution in a particular region of interest, a Chimera mesh allows one to use a fine foreground mesh that moves as the region of interest while elsewhere (e.g. in the background) the mesh will be coarser, thus optimizing computational effort and time. 

The classification of cells distributed over background and foreground partitions follows the clustering proposed by Sharma et al. in \cite{sharma2021overset}. Specifically, any cell can be grouped in one of this subsets:
\begin{itemize}
	\item \textit{Field cells}: cells over which all involved differential operators of the PDE are discretized. They are also called internal cells. The degrees of freedom related to this class of cells correspond to the degrees of freedom of the numerical discretization of the continuous differential equation.
	\item \textit{Fringe cells}: cells devoted to the communication among the different mesh blocks. They are located at the external boundary of foreground partitions and at the boundary of the hole in the background mesh.
	\item \textit{Hole cells}: inactive cells of background. For these cells, no local discretization of differential operators is performed and the degrees of freedom linked to them are not included in the numerical discretization of the continuous PDE.  
\end{itemize}
When the domain evolves in time, also the cell clustering is dynamic. Consequently, the same cell could belong to any of the three classes over time but always to one and only one of them at fixed time. For a given overset configuration, let $N_c$ be the number of active cells (i.e. internal or fringe). We denote with $\TOmega$ the union of all $N_c$ active cells in the background and foreground partitions, labeled with $\Tbg$ and $\Tfg$, respectively. It follows that $\TOmega$ is the tessellation employed for discretizing the computational domain $\Omega$. We remark that one has actually an explicit time dependency in the definition of $\TOmega$, namely $\TOmega(t)$, which is omitted in the sequel to make notation easier.

\subsection{Space and time discretization} \label{subsec.spacetime}
From now on, we explicitly denote the dependency of the computational domain on time $t$ as $\Omega(t)$, for $t \in [0, t_f] \subset \R^+$, with $t_f$ being a prescribed finite time. Let
\begin{equation}
	\Omegabg(t) = \bigcup_{\cell{i} \in \Tbg} \cell{i} \quad \text{and} \quad \Omegafg(t) = \bigcup_{\cell{i} \in \Tfg} \cell{i}
\end{equation}
be subdomains in $\Omega(t)$ defined by the union of active cells $\cell{i}$ in the background $\Tbg$ and foreground $\Tfg$ partitions, respectively. Even though no evolution is prescribed for the computational domain $\Omega$, if a deformation is imposed to the foreground subdomain $\Omegafg(t)$, the evolution affects also the background subdomain $\Omegabg(t)$ (because of the evolution of the hole) and the global domain $\Omega$, since it holds that $\Omegabg(t) \cup \Omegafg(t) = \TOmega(t) \subseteq \Omega(t)$ at any time $t$. In particular, if physical internal boundaries are affected by a displacement (e.g. solid bodies moving in a fluid), this implies that also the computational boundaries $\partial \Omega(t) = \Gamma_b(t) \cup \Gamma_c$ are evolving. Otherwise, if internal boundaries are neither evolving nor present but a displacement is prescribed to the foreground subdomain $\Omegafg(t)$, the computational boundaries fulfill the relation $\partial \Omega(t) \equiv \partial \Omega(0) = \Gamma_c$ for any $t \in [0, t_f]$. 

In order to properly introduce the numerical scheme, the time set $[0, t_f]$ is split in intervals $[t^n, t^{n+1}]$, with $n = 1, \ldots, N_t$, such that $t^0 = 0$, $t^{N_t} = t_f$ and 
\begin{equation} \label{eq.time}
	t^{n+1} = t^n + \dt, 
\end{equation}
where the time step size $\dt$ is computed at any time iteration from $t^n$ to $t^{n+1}$ in order to ensure stability of the numerical method, as defined by \eqref{eq.CFL} in Section \ref{sec.numerical_method}. Any generic variable $z(t)$ depending on time and evaluated at the discrete time $t^n$ is denoted as $z^n$. 

At fixed time instance $t^n$, the overset configuration $\TOmega^n$ (of cardinality $N_c$) discretizes the computational space $\Omega^n$ through active cells $\omega_i \in \TOmega^n$ of area $| \omega_i |$, characteristic size $h_i = \sqrt{| \omega_i |}$, boundary $\partial \omega_i = \bigcup_{k = 1}^4 \Gamma_{ik}$ and center of mass $\xx_i$ given by
\begin{equation} \label{eq.centermass}
	\xx_i = \frac{1}{| \cell{i} |} \int_{\cell{i}} \xx \, \diff \Omega,
\end{equation}
with $\xx=(x,y)$ being the generic position vector over the domain $\Omega$. 

\subsection{Stencil of active cells} \label{subsec.stencil}
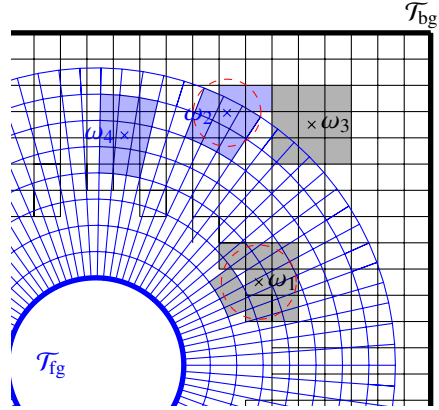
\begin{figure}
	\begin{center}
		\begin{tikzpicture}
			\begin{axis}[xmin=0,xmax=8,ymin=0,ymax=8,
				hide axis,
				axis equal,
				enlargelimits={abs=15pt},
				xlabel={x},
				ylabel={y},
				]
				
				\addplot [mark=none, black] table {overset_active_bg.dat};
				
				\addplot [mark=none, fill=black, opacity=0.3]  table {stencil_bg.dat};
				\addplot [mark = x, black] table {cell_bg.dat};
				\coordinate (Gamma3) at (axis cs:4.9,5.5);
				\draw (Gamma3) node[anchor=west]{$\cell{3}$};

				\addplot [mark=none, blue]  table {overset_active_fg.dat};
				
				\addplot [mark=none, fill=blue, opacity=0.3]  table {stencil_fg.dat};
				\addplot [mark = x, blue] table {cell_fg.dat};
				\coordinate (Gamma4) at (axis cs:0.65,5.25);
				\draw (Gamma4) node[anchor=east, text=blue]{$\cell{4}$};
				
				
				\addplot [mark=none, fill=black, opacity=0.3]  table {fringe_stencil_bg.dat};
				\addplot [mark = x, black] table {fringe_cell_bg.dat};
				\addplot [mark = none, red, dashed] table {fringe_stencil_circle_bg.dat};
				\coordinate (Gamma1) at (axis cs:3.7,1.9);
				\draw (Gamma1) node[anchor=west]{$\cell{1}$};
				
				\addplot [mark=none, fill=blue, opacity=0.3]  table {fringe_stencil_fg.dat};
				\addplot [mark = x, blue] table {fringe_cell_fg.dat};
				\addplot [mark = none, red, dashed] table {fringe_stencil_circle_fg.dat};
				\coordinate (Gamma2) at (axis cs:2.9,5.7);
				\draw (Gamma2) node[anchor=east,text=blue]{$\cell{2}$};
				
				\addplot [mark=none, black, line width=2pt]  table {Gamma_c.dat};
				
				\addplot [mark=none, blue, line width=2pt]  table {Gamma_b.dat};
				
				\coordinate (Gammac) at (axis cs:8,8);
				\draw (Gammac) node[anchor=east]{$\Tbg$};
				
				\coordinate (Gammab) at (axis cs:-1,0);
				\draw (Gammab) node[blue]{$\Tfg$};

			\end{axis}
		\end{tikzpicture}
	\end{center}
	\caption{Minimal stencils on both background and foreground cells. Background stencils are highlighted in black and foreground stencils in blue. For fringe cells $\cell{1} \in \Tbg$ and $\cell{2} \in \Tfg$, the relative stencils $\Ste{1}$ and $\Ste{2}$ are composed of cells belonging to both partitions. In particular, cells in the other partition have their center of mass internal to the stencil-circles (in dashed red lines for both stencils). For internal cells $\cell{3} \in \Tbg$ and $\cell{4} \in \Tfg$, the relative stencils $\Ste{3}$ and $\Ste{4}$ are composed of all cells sharing at least one vertex with $\cell{3}$ and $\cell{4}$, respectively.}
	\label{fig.stencils}
\end{figure}

Regardless of whether it belongs to a background or foreground block, any active cell $\cell{i}$ has a stencil $\Ste{i}$ of neighboring cells. The union of all cells in the stencil defines a local region of the space over which the numerical approximation of all involved integral or differential operators as well as interpolation or extrapolation processes are performed. The size of the stencil has an impact on the computational cost of these operations: the larger the stencil, the greater the computational effort to be made. Thus, we limit any stencil to be minimal, according to the following definitions.
\begin{definition}[\textbf{Stencil-circle}] \label{def.stencil-circle}
	A stencil-circle $\mathcal{C}_i$ of a stencil $\Ste{i}$ is a circle whose circumference is centered at $\xx_i$ with radius equal to the double of the maximum distance among the center of mass $\xx_i$ and any vertex of cell $\cell{i}$.  
\end{definition}
\begin{definition}[\textbf{Minimal stencil}] \label{def.minmal_stencil}
	A stencil $\Ste{i}$ is said to be minimal when it is composed of all cells $\cell{j}$, with $j \neq i$, fulfilling one of the two following options:
	\begin{itemize}
		\item cells $\cell{j}$ and $\cell{i}$ share at least one vertex;
		\item center of mass $\xx_j$ is internal to the stencil-circle $\mathcal{C}_i$.
	\end{itemize} 
\end{definition} 
If cell $\cell{i} \in \mathcal{T}_{\star}$ ($\star = \text{bg}, \text{fg}$) is internal, the minimal stencil $\Ste{i}$ is composed of all cells $\cell{j} \in \mathcal{T}_{\star}$ sharing at least one vertex with $\cell{i}$. As originally proposed in \cite{bergmann2022second}, for a fringe cell $\cell{i}$ in partition $\Tbg$ ($\Tfg$), the cells of stencil $\Ste{i}$ of the same partition are the ones sharing at least one vertex with $\cell{i}$. The remaining cells in the other partition $\Tfg$ ($\Tbg$) are chosen among the ones whose centers of mass are internal to the stencil-circle $\mathcal{C}_i$ but are not covered by any other cell in $\Tbg$ ($\Tfg$) of the stencil itself. In Figure \ref{fig.stencils} we sketch some examples of possible fringe and internal stencils on background and foreground meshes. 

\section{Governing equations} \label{sec.math}
Let $\uu(\xx,t) : \Omega(t) \times [0, t_f] \rightarrow \R^2$ and $p(\xx,t) : \Omega(t) \times [0, t_f] \rightarrow \R$ be the fluid velocity field and pressure of an incompressible fluid, respectively. They fulfill the Navier-Stokes equations
\begin{subequations} \label{eq.INS_phys}
	\begin{align}
		\nabla \cdot \uu &= 0 &&\text{ in } \Omega(t) \times [0, t_f], \label{eq.INS_Phys_cont}\\
		\frac{\partial \uu}{\partial t} + \nabla \cdot (\uu \otimes \uu) - \nu \Delta \uu + \nabla p &= \mathbf{0} &&\text{ in } \Omega(t) \times [0, t_f], \label{eq.INS_Phys_mom}\\
		\uu(\xx, 0) &= \uu_0(\xx) &&\text{ in } \Omega(0) \times \{ 0 \}, \label{eq.INS_Phys_IC}
	\end{align}
\end{subequations}
completed by proper boundary conditions on $\partial \Omega(t) \times [0, t_f]$. Equations \eqref{eq.INS_Phys_cont} and \eqref{eq.INS_Phys_mom} represent the balance of mass and momentum, respectively. The kinematic viscosity coefficient $\nu = \mu / \rho$ is given by the ratio between dynamic viscosity $\mu$ and density $\rho$, which is assumed to be constant. Consequently, the kinematic viscosity is also constant, since the dynamic viscosity is a physical property of the fluid. The initial conditions are imposed in \eqref{eq.INS_Phys_IC} through the function $\uu_0=\uu(\xx,0)$.

The space discretization $\TOmega(t) \subseteq \Omega(t)$ is supposed to evolve at least in its subset $\Omegafg(t)$. For this reason, the equation for the frame motion is governed by a Cauchy problem for the position vector $\xx(t) \in \Omegafg(t)$, that is
\begin{subequations} \label{eq.motion}
	\begin{align}
		\frac{\diff \xx}{\diff t} &= \w && \text{in } [0, t_f], \quad \xx \in \Omegafg, \\
		\xx(0) &= \xx_0 && t = 0,
	\end{align}
\end{subequations} 
where $\xx_0$ denotes the initial foreground configuration and $\w = \w(\xx, t)$ is the foreground mesh evolution velocity that is assumed to depend on $\xx$ and $t$. 

In order to identify the terms responsible of a possible stiffness of physical model, it is more convenient to study the dimensionless form of system \eqref{eq.INS_phys} according to \cite{boscheri2021high, boscheri2023all}. In fact, phenomena described by this set of equations can be characterized by different time scales. In such cases, it may become difficult to numerically treat the solution. An example is given by viscous-dominated flows leading to a severe CFL-type stability constraint which can be regarded as a stiffness of the problem. Let $L_0$, $T_0$ and $U_0 = L_0 / T_0$ be the characteristic length, time and velocity of the phenomenon under consideration, respectively. Through them, it is possible to rescale all physical quantities involved in \eqref{eq.INS_phys} and \eqref{eq.motion} as 
\begin{equation} \label{eq.dimensionless}
	\tilde{\xx} = \frac{\xx}{L_0}, \quad \tilde{t} = \frac{t}{T_0}, \quad \tilde{\uu} = \frac{\uu}{U_0} \quad \tilde{p} = \frac{p L_0}{\nu U_0}, \quad \tilde{\w} = \frac{\w}{U_0}.
\end{equation}
The dimensionless variables \eqref{eq.dimensionless} allow equations \eqref{eq.INS_phys} and \eqref{eq.motion} to be reformulated as
\begin{subequations} \label{eq.INS_dim}
	\begin{align}
		\nabla \cdot \uu &= 0 &&\text{ in } \Omega(t) \times [0, t_f], \label{eq.INS_dim_cont}\\
		\frac{\partial \uu}{\partial t} + \nabla \cdot (\uu \otimes \uu) - \frac{1}{\Rey} \Delta \uu + \nabla p &= \mathbf{0} &&\text{ in } \Omega(t) \times [0, t_f], \label{eq.INS_dim_mom}\\
		\frac{\diff \xx}{\diff t} &= \w && \text{in } [0, t_f], \quad \xx \in \Omegafg, \label{eq.INS_dim_motion}
	\end{align}
\end{subequations}
where the tilde symbol has been removed over all dimensionless quantities for easing the readability. The dimensionless number $\Rey = U_0 L_0 / \nu$ is the Reynolds number and it represents the ratio between inertial and viscous forces in the studied phenomenon. In the stiff limit of the model, i.e. for $\Rey \to 0$, the governing equations \eqref{eq.INS_dim} reduce to the Stokes equations at first order leading terms \cite{boscheri2023new}:
\begin{equation}
	\label{eqn.Stokes}
	\nabla p - \nu \Delta \uu = \mathbf{0}.
\end{equation}

For the sequel, we rewrite system \eqref{eq.INS_dim} in a compact form. Let 
\begin{equation} \label{eq.INS_strong_tuple}
	\QQ = ( 0, \uu ), \quad \F(\QQ, \nabla \QQ) = \left( \uu, \uu \otimes \uu - \frac{1}{\Rey}\nabla \uu + p \I \right),
\end{equation}
be the tuple of unknown variables and nonlinear flux, respectively. In \eqref{eq.INS_strong_tuple}, matrix $\I$ denotes the identity. Consequently, the strong formulation of the incompressible Navier-Stokes equations reads
\begin{subequations} \label{eq.INS_strong}
	\begin{align}
		\frac{\partial \QQ}{\partial t} + \nabla \cdot \F(\QQ, \nabla \QQ) &= \boldsymbol{0}, &&\xx \in \Omega(t), \quad t \in [0, t_f], \quad \QQ \in \Omega_{\QQ}, \label{eq.INS_strong_INS} \\
		\frac{\diff \xx}{\diff t} &= \w, && \xx \in \Omegafg(t), \quad t \in [0, t_f],  \label{eq.INS_strong_mot}
	\end{align}
\end{subequations}
with the vector of evolutionary variables $\QQ$ in \eqref{eq.INS_strong_INS} defined in the space $\Omega_{\QQ} \subset \R^\gamma$ of admissible states. 

\subsection{Frame evolution integration}
We introduce an integral version of equation \eqref{eq.INS_dim_mom} allowing to account for the evolution of frame $\xx \in \Omegafg$ directly inside the nonlinear flux term. Since the deformation velocity $\w$ only refers to the space coordinates in the foreground subdomain $\Omegafg$, let us extend it to the whole domain as
\begin{equation} \label{eq.w_extension}
	\tilde{\w} = \left\{ \begin{matrix}
		\w, & \text{in } \Omegafg \\
		\mathbf{0}, & \text{in } \Omegabg
	\end{matrix} \right. .
\end{equation}
By abuse of notation, the extended velocity in \eqref{eq.w_extension} is denoted without the tilde symbol for ease of reading. From now on, the velocity $\w$ refers to its extension \eqref{eq.w_extension}.

Let $\omega(t) \subseteq \Omega(t)$ be a generic and evolving in time control volume in $\Omega(t)$. Because of the Reynolds transport theorem, the integration over the control volume $\omega(t)$ of $\partial_t \uu$ in the left hand side of \eqref{eq.INS_dim_mom} is
\begin{equation} \label{eq.RTT}
	\int_{\omega(t)} \frac{\partial \uu}{\partial t} \, \diff \xx = \frac{\diff}{\diff t} \int_{\omega(t)} \uu \, \diff \xx - \int_{\partial \omega(t)} \uu \, ( \dot{\xx} \cdot \nn ) \, \diff \Gamma = \frac{\diff}{\diff t} \int_{\omega(t)} \uu \, \diff \xx - \int_{\partial \omega(t)} [\uu \otimes \w] \nn \, \diff \Gamma ,
\end{equation}
where $\nn$ is the outer unit normal to the space-time boundary $\partial \omega(t)$. It is now possible to introduce the integral version of the momentum equation \eqref{eq.INS_dim_mom} for a generic control volume $\omega(t)$ as
\begin{equation} \label{eq.INS_int}
	\frac{\diff}{\diff t} \int_{\omega(t)} \uu \, \diff \xx + \int_{\partial \omega(t)} \left[ \F^{\w}(\uu) - \frac{1}{\Rey} \nabla \uu + p \I \right] \nn \, \diff \Gamma = 0,
\end{equation}
where the associated nonlinear flux reads
\begin{equation} \label{eq.int_flux}
	\F^{\w}(\uu) = \uu \otimes (\uu - \w).
\end{equation}
We remark that system \eqref{eq.INS_int} is still coupled with the motion equation \eqref{eq.INS_strong_mot} restricted to the control volume $\omega(t)$. In our approach, the mesh velocity can be arbitrarily chosen, hence dealing with an Arbitrary-Lagrangian-Eulerian (ALE) description of the continuum.

\section{Numerical method} \label{sec.numerical_method}
This section is devoted to the numerical approach for integrating the Navier-Stokes system \eqref{eq.INS_strong}. We start by splitting the system in two sub-systems with respect to the scaling parameter defined by the Reynolds number. This allows to express a semi-discrete scheme in time. Successively, through a finite volume method, the system is fully discretized in space as well. The arising spatial discretizations are explained as well as the treatment of information on fringe cells of the overset configuration. Finally, we introduce the high-order IMEX approach on moving Chimera meshes. 

The numerical method belongs to the category of Arbitrary-Eulerian-Lagrangian (ALE) methods since the mesh velocity $\w$ in \eqref{eq.motion} can be arbitrarily chosen. The fully discrete scheme is second order accurate in both space and time.

\subsection{Flux splitting}
In order to address the challenges posed by multi-scale problems, as discussed in the previous section, as well as the movement of the mesh or the evolution of the computational domain, an efficient and accurate numerical method needs to be developed. To achieve this, a flux splitting approach is adopted. Consequently, the governing equations are separated according to different scales, namely fast and slow, with respect to the conservative flux terms. This strategy has been previously employed in the literature to handle both incompressible flow models \cite{casulli1990semi, boscheri2023all, toro2022flux, vater2018semi, boscheri2023new} and compressible fluids \cite{BOSCHERI2022127457,boscheri2021high, toro2012flux, dumbser2016conservative}. The core idea behind this technique is to partition the system of governing equations into two sub-systems: one that contains fluxes dependent on the scaling parameter (i.e. the Reynolds number), and the other one that does not. By doing so, this procedure effectively distinguishes the terms that will be explicitly discretized in time from those that require an implicit discretization due to their dependence on the Reynolds number. The flux splitting for system \eqref{eq.INS_strong_INS} reads
\begin{equation} \label{eq.split_strong}
	\frac{\partial \QQ}{\partial t} + \nabla \cdot \EX{\F}(\QQ) + \nabla \cdot \IM{\F}(\QQ, \nabla \QQ) = \mathbf{0}.
\end{equation}
Subscripts 'E' and 'I' refer to the explicit and the implicit terms, respectively. The split form \eqref{eq.split_strong} naturally induces to a partitioned system \cite{hofer1976partially} defined by a convective and a pressure-viscosity sub-system that write
\begin{subequations} \label{eq.EX_IM_strong}
	\begin{align}
		\frac{\partial \QQ}{\partial t} + \nabla \cdot \EX{\F}(\QQ) &= \mathbf{0}, & \text{(convective)}  \label{eq.EX_IM_conv_strong} \\
		\frac{\partial \QQ}{\partial t} + \nabla \cdot \IM{\F}(\QQ, \nabla \QQ) &= \mathbf{0}, & \text{(pressure-viscosity)} \label{eq.EX_IM_pres_strong}
	\end{align}
\end{subequations} 
with the explicit and implicit fluxes given by
\begin{equation} \label{eq.EX_IM_fluxes_strong}
	\EX{\F}(\QQ) = ( \mathbf{0}, \uu \otimes \uu ), \quad \IM{\F}(\QQ, \nabla \QQ) = \left( \uu, -\frac{1}{\Rey} \nabla \uu + p \I  \right).
\end{equation}

In order to properly define the time step size $\dt$, the explicit (convective) sub-system \eqref{eq.EX_IM_conv_strong} is analyzed. The eigenvalues in normal direction $\nn=(n_x,n_y)$ for sub-system \eqref{eq.EX_IM_conv_strong} are
\begin{equation} \label{eq.eigenvalues}
	\lambda_1 = 0, \quad \lambda_{2,3} = 2 ( \uu - \w ) \cdot \nn.
\end{equation}
For a computational mesh $\TOmega^n$ of characteristic size $h$ at time instance $t^n$, the time step size $\dt = t^{n+1} - t^n$ fulfills a classical CFL stability condition with respect to the maximum convective eigenvalue, thus
\begin{equation} \label{eq.CFL}
	\dt \leq \CFL \min_{\TOmega^n} \frac{h}{\max | \lambda |}.
\end{equation} 
We remark that the time step size $\dt$ is computed at each time iteration not only because the maximum convective eigenvalue evolves with the fluid velocity but also because the characteristic mesh size $h$ could evolve due to the mesh deformation velocity $\w$. Condition \eqref{eq.CFL} results in a less stringent stability requirement, especially in the asymptotic regime given by $\Rey \rightarrow 0$, compared to a fully explicit time discretization. This implies that the time step becomes independent of the stiffness parameter, allowing simulations of flows characterized by different Reynolds numbers to be executed with equal computational efficiency. Furthermore, because the terms related to the stiffness are treated implicitly, the numerical dissipation is solely proportional to the fluid speed. This feature makes the numerical scheme particularly well-suited for applications in the asymptotic limit (zero-relaxation) when the Reynolds number vanishes, that is the Stokes model \eqref{eqn.Stokes}. It is worth noting that problems of this nature could not be simulated using a purely explicit scheme \cite{guillard2004behavior, guillard1999behaviour, dellacherie2010analysis}.

\subsection{Finite volume space}
A finite volume scheme is used to discretize the explicit terms of system \eqref{eq.EX_IM_conv_strong} and also for solving the implicit sub-system \eqref{eq.EX_IM_pres_strong}, hence storing the solution within the control volume $\cell{i} \in \TOmega^n$ at any time instance $t^n$. In particular, the vector of conserved variables $\QQ$ is represented as cell averages referring to each cell element $\cell{i}$ as
\begin{equation} \label{eq.cell_average}
	\QQ_i^n := \frac{1}{| \cell{i} |} \int_{\cell{i}} \QQ(\xx, t^n) \, \diff \xx.
\end{equation}

For achieving second-order accuracy, the numerical solution in \eqref{eq.cell_average} must undergo a reconstruction process. Let $E_i = \cell{i} \cup \bigcup_{ \cell{j} \in \Ste{i} } \cell{j}$ be the subdomain defined by all cells in the stencil $\Ste{i}$ centered on cell $\cell{i}$ of characteristic size $h_i$. For a given function $\Phi \in \mathcal{C}^2(E_i)$, whose knowledge is restricted to the cell centers in $E_i$ through \eqref{eq.cell_average}, its quadratic polynomial reconstruction is denoted by $\Pol{i}(\Phi)$ and it is expressed using the Taylor polynomial space function
\begin{eqnarray} \label{eq.taylor}
	\mathcal{P}_2(E_i)&=&\text{span} \left\{ 1, \frac{x - x_i}{h_i}, \frac{y - y_i}{h_i}, \frac{(x - x_i)(y - y_i)}{h_i^2}, \frac{(x - x_i)^2}{2 h_i^2}, \frac{(y - y_i)^2}{2 h_i^2} \right\} \nonumber \\
	&=& \{z_{i0}^P,z_{i1}^P,z_{i2}^P,z_{i3}^P,z_{i4}^P,z_{i5}^P\}:=\zz^P_i(\xx),
\end{eqnarray}
that constitutes the set of basis functions $\zz^P_i(\xx)$. The final polynomial reconstruction is then obtained as an expansion of the form
\begin{equation} \label{eq.polyrec}
	\Pol{i}(\Phi) = (\zz^{P}_i)^\top \, \hat{\boldsymbol{\Phi}}_{i}.
\end{equation}
The unknown polynomial coefficients $\hat{\boldsymbol{\Phi}}_{i}=\{\hat{\Phi}_{i0},\hat{\Phi}_{i1},\hat{\Phi}_{i2},\hat{\Phi}_{i3},\hat{\Phi}_{i4},\hat{\Phi}_{i5}\}$ related to reconstruction $\Pol{i}(\Phi)$ are found by imposing that the polynomial exactly coincides with the function at the cell center $\xx_i$ (i.e., $\Pol{i} (\Phi(\xx_i)) \equiv \Phi_i$) and in the mean-square sense over all other cell centers $\cell{j} \in \Ste{i}$ (namely, $\Pol{i} (\Phi(\xx_j)) = \Phi_j $ for any $\cell{j} \in \Ste{i}$). This implies that the polynomial coefficient of the first component of the basis, corresponds to the polynomial evaluation at cell center $\xx_i$, i.e. $\hat{\Phi}_{i0}=\Phi_i$. Further explanations on this reconstruction is provided in Section \ref{sec.fringe} for the compact transmission of information from one partition to another in the overset configuration along the fringe cells.

\begin{remark}
	The polynomial space function $\mathcal{P}_2(E_i)$ in \eqref{eq.taylor} properly works for solutions that are continuous at least up to the second order derivative. It is employed in this work by exploiting the $\mathcal{C}^\infty$-continuity of the solution for problem \eqref{eq.INS_dim}.
\end{remark}

\subsection{First order semi-implicit scheme} \label{sec.IMEX1}
We start introducing the first order semi-implicit IMEX scheme in time. This gives the opportunity to explain the fundamental steps of integration in time of equation \eqref{eq.split_strong} over one time step. In the first order scheme, the current time is $t^n$ and the later time is the next time level $t^{n+1} = t^n + \dt$. 
The first order semi-implicit IMEX scheme for the splitting formulation \eqref{eq.split_strong} is written as
\begin{equation} \label{eq.IMEX1}
	\frac{\partial \QQ}{\partial t} + \nabla \cdot \EX{\F}(\QQ^n) + \nabla \cdot \IM{\F}(\QQ^{n+1}, \nabla \QQ^{n+1}) = \mathbf{0} \quad \text{in } \Omega(t) \times [t^n, t^{n+1}],
\end{equation}
which explicitly yield
\begin{subequations} \label{eq.INS_IMEX1}
	\begin{align}
		\nabla \cdot \uu^{n+1} &= 0, \label{eq.INS_IMEX1_cont}\\
			\frac{\partial \uu}{\partial t} + \nabla \cdot (\uu^n \otimes \uu^n) - \frac{1}{\Rey} \Delta \uu^{n+1} + \nabla p^{n+1} &= \mathbf{0}. \label{eq.INS_IMEX1_mom}
	\end{align}
\end{subequations}
Inspired by \cite{bergmann2022ader}, the above system is solved by a fractional-step method in the paradigm of the prediction-projection approach. Firstly, in the prediction step, an intermediate velocity $\ustar$ is computed by not considering the divergence-free constraint. Successively, a correction step is carried out in order to project the eventually nonphysical intermediate velocity onto a solenoidal space. This allows to ensure constraint \eqref{eq.INS_IMEX1_cont} to be respected at the semi-discrete level. 
The sequential steps of this approach are the prediction step
\begin{equation} \label{eq.prediction}
	\frac{\partial \ustar}{\partial t} - \nabla \cdot (\uu^n \otimes \uu^n) - \frac{1}{\Rey} \Delta \ustar + \nabla p^n = \mathbf{0} \quad \text{in } \Omega(t) \times [t^n, t^{n+1}],
\end{equation}
and the correction step
\begin{subequations} \label{eq.correction}
	\begin{align}
		\nabla \cdot \uu^{n+1} &= 0 & \text{in } \Omega^{n+1},  \label{eq.correction_constraint} \\
		\frac{\uu^{n+1} - \ustar}{\dt} + \nabla p^{n+1} - \nabla p^n &= \mathbf{0} & \text{in } \Omega^{n+1}. \label{eq.correction_eq}
	\end{align}
\end{subequations}
In the prediction step, the nonlinear flux refers to the current time $t^n$, as suggested by the flux splitting approach in \eqref{eq.IMEX1}. The pressure $p^{n+1}$ at next time appears in the correction step. However, pressure $p^n$ at time $t^n$ is taken into account as force term in \eqref{eq.prediction} for giving a physical information on the pressure field, as proposed in \cite{tavelli2014staggered}. After integration over a field cell $\cell{i}(t)$ and due to the Reynolds transport theorem \eqref{eq.RTT}, the prediction step \eqref{eq.prediction} is solved in the sense of finite volume as
\begin{equation} \label{eq.ustar}
	\int_{\cell{i}^{n+1}} \ustar \, \diff \xx - \frac{\dt}{\Rey} \int_{\partial \cell{i}^{n+1}} [ \nabla \ustar ] \, \nn \, \diff \Gamma = \f_{\uu}^n - \dt \int_{\partial \cell{i}^n} p^n \, \nn \, \diff \Gamma, \qquad \f_{\uu}^n = \int_{\cell{i}^n} \uu^n \, \diff \xx - \dt \int_{\partial \cell{i}^n} [\F^{\w^n}(\uu^n)] \, \nn \, \diff \Gamma,
\end{equation}
where the nonlinear flux $\F^{\w^n}$ considers the mesh velocity $\w^n$ at time $t^n$ through the definition \eqref{eq.int_flux}. The discrete convection term is then referred to as $\f_{\uu}^n$.
The resolution of problem \eqref{eq.ustar} needs the knowledge of the space evolution from $t^n$ to $t^{n+1}$ according to the motion equation \eqref{eq.INS_dim_motion}. This Cauchy problem is implicitly solved as
\begin{equation} \label{eq.space_problem}
	\frac{\xx^{n+1} - \xx^n}{\dt} = \w^{n+1}.
\end{equation} 
The above equation could be nonlinear since $\w^{n+1}=\w^{n+1}(\xx^{n+1},t^{n+1})$, hence we resort to a Newton method for linearizing and solving it. The intermediate velocity $\ustar$ solving \eqref{eq.ustar} contains the computation of both explicit convection and implicit viscous contribution as well as the information on the evolution of the frame, obtained from \eqref{eq.space_problem}. 
Once the provisional velocity $\uu^*$ is known, it is possible to treat the correction step \eqref{eq.correction}.
By plugging the second equation \eqref{eq.correction_eq} in the energy equation \eqref{eq.correction_constraint}, we obtain an elliptic equation for the resolution of the pressure field $p^{n+1}$:
\begin{equation} \label{eq.pressure}
	\dt \, \Delta p^{n+1} \,  = \nabla \cdot \ustar + \dt \, \Delta p^n.  
\end{equation}
The elliptic problem \eqref{eq.pressure} is solved relying on a finite volume approach, hence obtaining
\begin{equation} \label{eq.pressure_int}
	\dt \, \int_{\cell{i}^{n+1}}  \Delta p^{n+1} \, \diff \xx = \int_{\partial \cell{i}^{n+1}} \ustar \cdot \nn \, \diff \Gamma + \dt \, \int_{\cell{i}^{n+1}} \Delta p^n \, \diff \xx.  
\end{equation}
Once the new pressure $p^{n+1}$ is known, the velocity field $\uu^{n+1}$ is updated by \eqref{eq.correction_eq}. This closes the semi-discretization of problem \eqref{eq.IMEX1}.

\begin{remark}[Divergence-free constraint]
	The presented semi-discrete scheme provides a divergence-free velocity field, i.e. $\nabla \cdot \uu^{n+1} = 0$. As a matter of fact, it holds that
	\begin{equation} \label{eq.helmholtz}
		\ustar = \uu^{n+1} + \dt \, \nabla ( p^n - p^{n+1} ).
	\end{equation} 
	In the above equation, the velocity $\uu^{n+1}$ is solenoidal by imposition of constraint \eqref{eq.correction_constraint} and the term $\dt \nabla ( p^n - p^{n+1} )$ is irrotational by definition. This means that, under usual regularity assumption of all variables at boundary $\partial \Omega$, the velocity $\ustar$ decomposes uniquely according to Helmholtz decomposition theorem. Thus, system \eqref{eq.correction} is a projection step of the intermediate velocity over a solenoidal space.
\end{remark}

\begin{remark}[Aribitrary Lagrangian-Eulerian approach]
	In the Arbitrary Lagrangian-Eulerian context, the usage of an intermediate velocity $\ustar$ from the presented fractional-step method permits to provide to both pressure $p^{n+1}$ (through \eqref{eq.pressure}) and velocity field $\uu^{n+1}$ (via \eqref{eq.correction}) the knowledge of frame evolution within a splitting scheme. 
\end{remark}

\begin{remark}[Fractional-step method]
	The introduced method is an incremental Chorin-Temam type fractional-step method \cite{chorin1968numerical, temam1969approximation, quarteroni1999analysis}, since it uses the notion of the previous time pressure $p^n$ to obtain an intermediate velocity $\ustar$ and then a current time pressure $p^{n+1}$ to project the velocity (possibly non-solenoidal and therefore nonphysical) to a divergence-free manifold.
\end{remark}

Now, let us list the most salient points of the new scheme at first order:
\begin{enumerate}
	\item Solve the motion equation \eqref{eq.INS_dim_motion} through an implicit scheme \eqref{eq.space_problem}. If the arising algebraic problem is nonlinear, linearize it using Newton method.
	\item Solve the prediction step \eqref{eq.ustar} for the intermediate velocity field $\ustar$.
	\item Through the divergence of the previously computed velocity $\ustar$, solve the elliptic equation for pressure $p^{n+1}$ given by \eqref{eq.pressure_int}.
	\item Finally, update the velocity field $\uu^{n+1}$ via \eqref{eq.correction_eq} through the intermediate velocity $\ustar$ and the gradient of the new pressure $p^{n+1}$. 
\end{enumerate}

\subsection{Spatial discretization of the motion equation} \label{sec.trajODE}
The trajectory equation \eqref{eq.INS_dim_motion} is defined at each vertex $v_k$ of the computational mesh of coordinates $\xx_{vk}$, thus the semi-discrete scheme \eqref{eq.space_problem} is spatially approximated as
	\begin{equation} \label{eq.spacetime_problem}
		\frac{\xx^{n+1}_{v_k} - \xx^n_{v_k}}{\dt} = \w^{n+1}_{v_k}.
	\end{equation}
If the mesh velocity depends on $\xx$, i.e. $\w^{n+1}_{v_k}=\w^{n+1}_{v_k}(\xx^{n+1}_{v_k},t^{n+1})$, the nonlinear algebraic equation is solved at the aid of a Newton method. Once the mesh velocity is computed at each vertex of the computational mesh, then the new mesh configuration is simply given by the new coordinates $\xx^{n+1}_{v_k}$. 

In this way, we maintain the all control volumes in the mesh defined by straight lines, thus keeping quadrilateral elements $\cell{i}$ in the mesh that indeed leads to a second order approximation of the geometry at any time. As such, the mesh velocity can be easily interpolated over the control volumes. In particular, the mesh velocity at the midpoint $\xx_{ij}$ of an edge $\Gamma_{ij}$ shared between cells $\cell{i}$ and $\cell{j}$ is given by
\begin{equation}
	\label{eqn.w_edge}
	\w_{ij} = \frac{1}{2} (\w_{v_1}+\w_{v_2}),
\end{equation}
where $v_1,v_2$ denote the nodes attached to edge $\Gamma_{ij}$, as depicted in Figure \ref{fig.diamond}.

\subsection{Spatial discretization of the Laplacian operator} \label{sec.laplacian}
\begin{figure}
	\centering
	\begin{tikzpicture}			
		\draw (0,0) coordinate (v1) node[below,xshift=.2cm] {\textcolor{red}{$v_1$}} -- (1.2,2) coordinate (v2) node[above,xshift=-.2cm] {\textcolor{red}{$v_2$}} -- (-1.5,2.1) -- (-1.9,.5) node[below] {$\omega_i$} -- (v1);
		\draw (0,0) -- (2,0) -- (3,2.5) node[right] {$\Omega_j$} -- (1.2,2) -- (0,0);
		
		\draw (.6,1) coordinate (P) node {$\bullet$};
		\draw[blue,thick] (-.55,1.15) coordinate (c1) node {$\circ$} node[left] {$\xx_i$} -- (1.3,1.125) coordinate (c2) node {$\circ$} node[right] {$\xx_j$};
		\draw[blue, ->, thick] (c1) -- (c2) node[above] {$\mathbf{c}_{ij}$};
		
		\draw[ultra thick,->] (P) -- (1.5,0.46) node[below] {$\boldsymbol{n}$};
		
		\draw[red,dashed,thick] (-.6,-1) -- (v1);
		\draw[red,thick,->] (v1) node[above, yshift=.05cm, xshift=-.01cm] {$\Gamma_{ij}$} -- (v2) node[right] {$\boldsymbol{\tau}_{ij}$};
		\draw[red,dashed,thick] (v2) -- (1.8,3) ;
		
		\draw[fill=gray, opacity=0.3] (c1) -- (c2) -- (1.2,-.5) -- (-1.4,-.55) node[left] {$\tilde{\mathcal{S}}_{v_1}$} -- (c1);
		
		\draw[fill=gray, opacity=0.3] (c1) -- (c2) -- (2.5,3) -- (0.3,2.7) node[left] {$\tilde{\mathcal{S}}_{v_2}$} -- (c1); 
		
	\end{tikzpicture}
	\caption{Sketch of two internal cells $\omega_i$ and $\omega_j$ sharing the edge $\Gamma_{ij}$. The filled cells in gray are the dual cells with respect to vertexes $v_1$ and $v_2$. Their vertexes are defined by the centers of mass of cells sharing the respective vertexes $v_1$ and $v_2$.}
	\label{fig.diamond}
\end{figure}
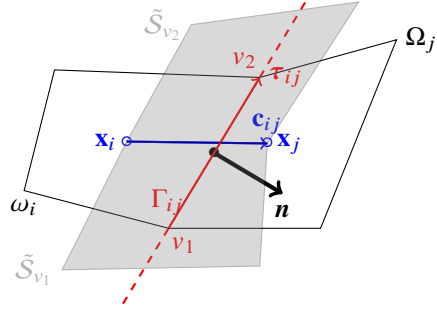

In both equations \eqref{eq.ustar} and \eqref{eq.pressure_int}, there is the necessity of discretizing the Laplacian operator applied to the solution $\Phi$, with $\Phi$ being the pressure field $p^{n+1}$ or one of the two components $(u^*, v^*)$ of the intermediate velocity $\ustar$. In the sense of finite volume, this means discretizing the quantity
\begin{equation} \label{eq.int_laplacian}
	\int_{\partial \cell{i}^{n+1}} \nabla \Phi \cdot \nn \, \diff \Gamma.
\end{equation} 
Let $\cell{i} \in \mathcal{T}_{\star}$, with $\star = \text{bg}, \text{fg}$, be an internal cell sharing with $\cell{j} \in \mathcal{T}_{\star}$ the edge $\Gamma_{ij} = \cell{i} \cap \cell{j}$. We denote by $\cc_{ij}$ the unit distance vector between the cell centers of the two neighboring cells, and with $\boldsymbol{\tau}_{ij}$ the unit normal vector tangent to $\Gamma_{ij}$. A possible discretization of the normal gradient on edge $\Gamma_{ij}$ exploits the second order accurate diamond formula \cite{bertolazzi2004cell, coudiere1999convergence}, that reads
\begin{equation} \label{eq.diamond}
	[\nabla \Phi \cdot \nn_{ij}]_{\Gamma_{ij}} \simeq \frac{1}{\cc_{ij} \cdot \nn_{ij}} \left( \frac{\Phi_j - \Phi_i}{| \xx_i - \xx_j |} - \boldsymbol{\tau}_{ij} \cdot \cc_{ij} \frac{\Phi_{v_2} - \Phi_{v_1}}{| \Gamma_{ij} |} \right), 
\end{equation} 
where $\Phi_{v_k}$, $k = 1,2$, is the evaluation of function $\Phi$ along vertexes $v_1$ and $v_2$, namely the extreme points of edge $\Gamma_{ij}$, as sketched in Figure \ref{fig.diamond}. Approximation \eqref{eq.diamond} is based on the fact that any vector (gradient included) can be expressed as a linear combination of a vector basis $\{ \cc_{ij}, \boldsymbol{\tau}_{ij} \}$. The reader is addressed to \cite{carlino2021ader} for the derivation of this approximation. Furthermore, when the cells are Cartesian (i.e., $\{ \cc_{ij}, \boldsymbol{\tau}_{ij} \}$ forms an orthonormal basis), the diamond formula turns into the second order centered finite difference operator. 

Since for any variable $\Phi$, its knowledge is confined to the cell center of a cell, an extrapolation is needed for expressing $\Phi_{v_k}$ in function of the values $\Phi$ at cell centers. For this reason, we use an extrapolation based on information deposed on all cells sharing vertex $v_k$, as originally proposed in \cite{bergmann2022ader}. Let $\tilde{\mathcal{S}}_{v_k}$ be the stencil of all cells sharing the vertex $v_k$ (as depicted by filled cells in Figure \ref{fig.diamond}). We denote with $\tilde{E}_{v_k} = \bigcup_{\cell{j} \in \tilde{\mathcal{S}}_{v_k}} \cell{j}$. The edges linking the cell centers of $\tilde{E}_{v_k}$ define the dual cell $\tilde{\omega}_{v_k}$ with cell center $\tilde{\xx}_{v_k} = (\tilde{x}_{v_k}, \tilde{y}_{v_k})$ and characteristic size $\tilde{h}_{v_k} = \sqrt{| \tilde{\omega}_{v_k} |}$. We remark that the cell center $\tilde{\xx}_{v_k}$ of the dual cell does not necessarily coincide with the vertex position $\xx_{v_k}$ (this happens only for Cartesian square cells because $\xx_{v_k}$ is the real center of mass for $\tilde{\omega}_{v_k}$, thus $\tilde{\xx}_{v_k}=\xx_{v_k}$). We introduce the bilinear polynomial space as
\begin{equation} \label{eq.Q1}
	\mathcal{Q}_1(\tilde{E}_{v_k}) = \text{span} \left\{ 1, \frac{x - \tilde{x}_{v_k}}{h_{v_k}} , \frac{y - \tilde{y}_{v_k}}{h_{v_k}} , \frac{(x - \tilde{x}_{v_k}) (y - \tilde{y}_{v_k})}{h_{v_k}^2} \right\}:=\zz^Q_{v_k}(\xx),
\end{equation}
with $\zz^Q_{v_k}=\{z_{i0}^Q,z_{i1}^Q,z_{i2}^Q,z_{i3}^Q,z_{i4}^Q\}$ representing the set of bilinear basis functions. The polynomial representation of $\Phi_{v_k}$, indicated with $\Polb{v_k} (\Phi_{v_k})$, is given in terms of the the basis functions defined by the space \eqref{eq.Q1}, and it explicitly writes
\begin{equation} \label{eq.pol_Q1}
	\Polb{v_k} (\Phi_{v_k}) = (\zz^{Q}_{v_k})^\top \, \hat{\boldsymbol{\alpha}}_{v_k}.
\end{equation}
The polynomial coefficients $\hat{\boldsymbol{\alpha}}_{v_k}=\{\hat{\alpha}_{v_k0},\hat{\alpha}_{v_k1},\hat{\alpha}_{v_k2},\hat{\alpha}_{v_k3},\hat{\alpha}_{v_k4}\}$ are computed by imposing that the polynomial approximation coincides with the exact values of function $\Phi$ at the cell centers in $\tilde{E}_{iv_k}$, namely $\Polb{v_k}( \Phi(\xx_j)) \equiv \Phi(\xx_j)$ for any $\cell{j} \in \tilde{\mathcal{S}}_{v_k}$. This yields a linear system 
\begin{equation}
	\mathbf{A}_{v_k} \hat{\boldsymbol{\alpha}}_{v_k} = \boldsymbol{\Phi}_{v_k},
\end{equation}
with $\mathbf{A}_{v_k} \in \R^{4 \times 4}$ and $\boldsymbol{\Phi}_{v_k} \in \R^4$ the vector of values of function $\Phi$ evaluated at the cell centers of $\tilde{E}_{v_k}$. This permits to rewrite the diamond formula \eqref{eq.diamond} as
\begin{equation} \label{eq.diamond2}
	[\nabla \Phi \cdot \nn_{ij}]_{\Gamma_{ij}} \simeq \frac{1}{\cc_{ij} \cdot \nn_{ij}} \left( \frac{\Phi_j - \Phi_i}{| \xx_i - \xx_j |} - \boldsymbol{\tau}_{ij} \cdot \cc_{ij} \frac{(\zz^{Q}_{v_2})^\top \mathbf{A}_{v_2}^{-1} \boldsymbol{\Phi}_{v_2} - (\zz^{Q}_{v_1})^\top \mathbf{A}_{v_1}^{-1} \boldsymbol{\Phi}_{v_1}}{| \Gamma_{ij} |} \right),
\end{equation}
where all quantities related to $\Phi$ only refer to the cell centers. The diamond formula \eqref{eq.diamond2} considers both functional and geometrical information of the cells sharing an edge. For a stencil $\Ste{i}$ on a field cell $\cell{i}$, it is always possible to completely perform the following set separation: $\Ste{i} = \Ste{i}^+ \cup \Ste{i}^{\times}$, with $\Ste{i}^+ \cap \Ste{i}^{\times} = \emptyset$, where $\Ste{i}^+$ and $\Ste{i}^{\times}$ collect all cells in $\Ste{i}$ sharing at least one edge with $\cell{i}$ and only one vertex with $\cell{i}$, respectively. With this notation, the integrated Laplacian operator in \eqref{eq.int_laplacian} is approximated as
\begin{equation} \label{eq.laplacian_approx}
	\int_{\partial \cell{i}^{n+1}} \nabla \Phi \cdot \nn \, \diff \Gamma \simeq \sum_{\cell{j}^{n+1} \in \Ste{i}^+} \frac{| \Gamma_{ij} |}{\cc_{ij} \cdot \nn_{ij}} \left( \frac{\Phi_j - \Phi_i}{| \xx_i - \xx_j |} - \boldsymbol{\tau}_{ij} \cdot \cc_{ij} \frac{(\zz^{Q}_{v_2,j})^\top \mathbf{A}_{v_2,j}^{-1} \boldsymbol{\Phi}_{v_2,j} - (\zz^{Q}_{v_1,j})^\top \mathbf{A}_{v_1,j}^{-1} \boldsymbol{\Phi}_{v_1,j}}{| \Gamma_{ij} |} \right)=:\mathbb{K}_h(\Phi),
\end{equation}
where we have introduced the notation $\mathbb{K}_h(\Phi)$ to compactly address the numerical discretization of \eqref{eq.int_laplacian}.

\subsection{Spatial discretization of the convective flux} \label{sec.flux}
The convective flux needed to be discretized in \eqref{eq.ustar} is 
\begin{equation} \label{eq.conv_flux}
	\int_{\partial \cell{i}^n} \F^{\w^n}(\uu^n) \, \nn \, \diff \Gamma = \int_{\partial \cell{i}^n} [ \uu^n \otimes (\uu^n - \w^n) ] \, \nn \, \diff \Gamma,
\end{equation}
that is used in the definition of the explicit term $\f_{\uu}^n$ for obtaining the intermediate velocity $\ustar$. We choose to employ a simple Rusanov-type numerical flux function, hence obtaining
\begin{equation} \label{eq.ex_term}
	\f_{\uu}^n = \int_{\cell{i}^n} \uu^n \, \diff \xx - \dt \sum_{ \cell{j}^n \in \Ste{i}^+ } \int_{ \Gamma_{ij} } \mathcal{F}^{\w^n}( \uu_{ij}^{-,n}, \uu_{ij}^{+,n} , \nn_{ij} ) \, \diff \Gamma,
\end{equation}
where $\mathcal{F}^{\w^n}$ is the local flux approximation combining information from the left and right velocity states $\uu_{ij}^{-,n}$, $\uu_{ij}^{+,n}$, respectively, and the unit outward normal vector $\nn_{ij}$ with respect of the edge $\Gamma_{ij}$ as well as its deformation trough the edge velocity $\w^n_{ij}$. Also in this case, the cell $\cell{i}$ is supposed to be internal. In particular, the Rusanov-type numerical flux is defined as 
\begin{equation} \label{eq.rusanov}
	\mathcal{F}^{\w^n}( \uu_{ij}^{-,n}, \uu_{ij}^{+,n} , \nn_{ij} ) = \frac{1}{2} \left( \F^{\w^n}(\uu_{ij}^{-,n}) + \F^{\w^n}(\uu_{ij}^{+,n}) \right) \cdot \nn_{ij} - \frac{1}{2} | s_{\max} | \left( \uu_{ij}^{+,n} - \uu_{ij}^{-,n} \right).
\end{equation} 
The amount of numerical dissipation $|s_{\max}|$ is set equal to the maximum eigenvalue \eqref{eq.eigenvalues} of the convective sub-systems related to the left and right state $\uu_{ij}^{\pm,n}$. The boundary integral is numerically approximated by the midpoint rule which achieves second order of accuracy. Concerning the states, they are defined as the boundary extrapolated values along the midpoint of edge $\Gamma_{ij}$ computed through the polynomial reconstruction \eqref{eq.taylor} within cells $\cell{i}$ and $\cell{j}$ sharing edge $\Gamma_{ij}$, i.e.
\begin{equation} \label{eq.states}
	\uu_{ij}^{-,n} = \Pol{i} (\Phi) |_{ \Gamma_{ij} }, \quad \uu_{ij}^{+,n} = \Pol{j} (\Phi) |_{\Gamma_{ij}}.
\end{equation}
Finally, the mesh velocity $\w_{ij}^n$ at the edge midpoint is computed using the linear interpolation \eqref{eqn.w_edge}. The discrete convective operator is compactly addressed with $\mathbb{F}_h(\uu,\w)$:
\begin{equation} \label{eqn.convOp}
	\mathbb{F}_h(\uu,\w):= \sum_{ \cell{j}^n \in \Ste{i}^+ } \int_{ \Gamma_{ij} } \mathcal{F}^{\w^n}( \uu_{ij}^{-,n}, \uu_{ij}^{+,n} , \nn_{ij} ) \, \diff \Gamma.
\end{equation}

\subsection{Spatial discretization of the pressure gradient and the velocity divergence operators} \label{sec.gradp_divu}
In the right hand side of equation \eqref{eq.ustar} and \eqref{eq.pressure_int} we have the pressure gradient $\nabla p^n$ and the divergence of the intermediate velocity $\nabla \cdot \ustar$ to be discretized, respectively. As previously assumed, we still hypothesize cell $\cell{i}$ to be a field cell. In the finite volume sense, by means of Gauss theorem, they read
\begin{equation} \label{eq.gradp_divu}
	\int_{\partial \cell{i}^n} p^n \nn \, \diff \Gamma, \qquad \text{ and } \qquad \int_{\partial \cell{i}^{n+1}} \ustar \cdot \nn \, \diff \Gamma,
\end{equation} 
with the outward unit vector $\nn$ of the cell boundary $\partial \cell{i}^{\tau}$. Here, $\tau=\{n,n+1\}$ according to the time level of the integral argument in \eqref{eq.gradp_divu}. Let $\Phi$ be either the pressure $p^n$ or one component of velocity $\ustar=(u^*,v^*)$, so that the integrals in \eqref{eq.gradp_divu} turn into
\begin{equation} \label{eq.abstract_gradp_divu}
	\int_{\partial \cell{i}^{\tau}} \Phi \, \nn \, \diff \Gamma.
\end{equation}
 In this case, we still exploit the knowledge of internal stencils. As such, the discretization of \eqref{eq.abstract_gradp_divu} is given by
\begin{equation} \label{eq.disc_gradp_divu}
	\int_{\partial \cell{i}^{\tau}} \Phi \, \nn \, \diff \Gamma \simeq \sum_{\cell{j} \in \Ste{i}^+} | \Gamma_{ij} | \frac{1}{2} \left( \Pol{i} (\Phi)|_{\Gamma_{ij}} + \Pol{j} (\Phi)|_{\Gamma_{ij}} \right) \, \nn=:\mathbb{G}_h(\Phi),
\end{equation}
where we use a central flux function relying on the reconstructed data evaluated along the edge $\Gamma_{ij}$ at the midpoint. Here, the midpoint rule is again adopted to numerically approximate the boundary integral in \eqref{eq.abstract_gradp_divu}. Furthermore, the abbreviation $\mathbb{G}_h(\Phi)$ has been introduced to denote the gradient operator \eqref{eq.disc_gradp_divu}.

\subsection{Treatment of algebraic systems on fringe cells} \label{sec.fringe}
In  Section \ref{sec.laplacian}, \ref{sec.flux} and \ref{sec.gradp_divu}, we always suppose the cell over which the discretization is performed to be internal, i.e. belonging to the class of field cells. Let the general integral-differential problem $L[\Phi] = f$ be an abstract representation of either problem \eqref{eq.ustar} or \eqref{eq.pressure_int}, with $L$ a specific operator applied to the solution $\Phi$ and $f$ an abstract representation of the known explicit terms. The numerical approximation of the abstract problem on the whole discretized space $\TOmega$ is 
\begin{equation} \label{eq.alg}
	\mathbf{L} \hat{\boldsymbol{\Phi}} = \mathbf{f},
\end{equation}
with $\mathbf{L} \in \R^{N_c \times N_c}$ denoting the discretization of $L$, $\hat{\boldsymbol{\Phi}} \in \R^{N_c}$ representing the numerical solution at the cell centers of all $N_c$ active cells, and $\mathbf{f} \in \R^{N_c}$ being a suitable discretization of $f$. For instance, matrix $\mathbf{L}$ and vector $\mathbf{f}$ can be given by the approximation of the Laplacian operator in \eqref{eq.laplacian_approx} or the gradient operator in \eqref{eq.disc_gradp_divu} on internal cells. We are now interested in the discretization of the algebraic system \eqref{eq.alg} on fringe cells relying on the reconstruction polynomial $\Pol{i}$ given by \eqref{eq.polyrec}, which here exploits the information on the neighboring active internal cells on the other partition. Firstly, let us consider a cell $\cell{i}$ endowed of a stencil $\Ste{i}$ of cardinality $N_s$ with a function $\Phi \in \mathcal{C}^2(E_i)$, with $E_i = \cell{i} \cup \bigcup_{\cell{j} \in \Ste{i}} \cell{j}$. The reconstruction polynomial is then formally expressed as \eqref{eq.polyrec}, that is
\begin{equation}
	\label{eqn.polyfringe}
	\Pol{i}(\Phi) = (\zz^{P}_i)^\top \, \hat{\boldsymbol{\beta}}_{i},
\end{equation}
where $\hat{\boldsymbol{\beta}}_{i}=\{\hat{\beta}_{i0},\hat{\beta}_{i1},\hat{\beta}_{i2},\hat{\beta}_{i3},\hat{\beta}_{i4},\hat{\beta}_{i5}\}$ are the sought unknown expansion coefficients of the reconstruction of the fringe cell $\cell{i}$ and $\zz^{P}_i(\xx)$ are the Taylor basis defined by \eqref{eq.taylor}. Next, we impose that the polynomial exactly coincides with the value of the interpolated function at $\xx_i$, hence obtaining $\hat{\beta}_{i0}=\Phi_i$. We rewrite expansion \eqref{eqn.polyfringe} as
\begin{equation} \label{eq.pol_expansion2}
	\Pol{i}(\Phi) = \Phi_i + (\tilde{\zz}^P_i)^\top \tilde{\boldsymbol{\beta}}_i,
\end{equation}
with $\tilde{\zz}^P_i=\{z_{i1}^P,z_{i2}^P,z_{i3}^P,z_{i4}^P,z_{i5}^P\} \in \R^5$ collecting the basis components in \eqref{eq.taylor} of order greater or equal to 1 and $\tilde{\boldsymbol{\beta}}_i=\{\hat{\beta}_{i1},\hat{\beta}_{i2},\hat{\beta}_{i3},\hat{\beta}_{i4},\hat{\beta}_{i5} \} \in \R^5$ defined by the corresponding expansion coefficients. In order to determine $\tilde{\boldsymbol{\beta}}_i$, we impose the polynomial \eqref{eq.pol_expansion2} to coincide with the function $\Phi$ on the remaining cell centers of $E_i$ in the sense of mean-squares, i.e. $\Pol{i}(\Phi(\xx_j)) = \Phi_j$ for any $\cell{j} \in \Ste{i}$. This means
\begin{equation} \label{eq.mean_square}
	(\tilde{\zz}^P_i(\xx_j))^\top \tilde{\boldsymbol{\beta}}_i = \Phi_j - \Phi_i, \quad \forall \cell{j} \in \Ste{i}.
\end{equation}
By collecting rows of $(\tilde{\zz}^P_i(\xx_j))^\top$ in matrix $\tilde{\mathbf{Z}}_i^P \in \R^{N_s \times 5}$, problem \eqref{eq.mean_square} is equivalent to
\begin{equation} \label{eq.mean_square2}
	\tilde{\mathbf{Z}}_i^P \, \tilde{\boldsymbol{\beta}}_i = \boldsymbol{\Phi}_j - \boldsymbol{\Psi}_i,
\end{equation}
with $\boldsymbol{\Phi}_j$ and $\boldsymbol{\Psi}_i$ in $\R^{N_s}$ defined as 
\[
	\boldsymbol{\Phi}_j = \{\Phi_j\} \in \R^{N_s} \quad \forall \cell{j} \in \Ste{i} \qquad \text{and} \qquad \boldsymbol{\Psi}_i = \{\Phi_i\} \in \R^{N_s}.
\]
The solution of the linear system \eqref{eq.mean_square2} is
\begin{equation} \label{eq.pol_coeff}
	\tilde{\boldsymbol{\beta}}_i = \tilde{\mathbf{Z}}_i^{P,\dagger} ( \boldsymbol{\Phi}_j - \boldsymbol{\Psi}_i ),
\end{equation}
where $\tilde{\mathbf{Z}}_i^{P,\dagger} = ((\tilde{\mathbf{Z}}_i^{P})^\top \tilde{\mathbf{Z}}_i^{P})^{-1} \, (\tilde{\mathbf{Z}}_i^{P})^\top$ represents the pseudoinverse matrix of the mean-square problem \eqref{eq.mean_square2}.

Concerning the original algebraic problem \eqref{eq.alg}, let $\cell{k} \in \mathcal{T}_{\star}$, $\star = \text{bg}, \text{fg}$, be a fringe cell. Moreover, let $\cell{i} \in \tilde{\mathcal{T}}_{\star}$, an internal cell on the other partition $\tilde{\mathcal{T}}_{\star}$ whose cell center minimizes the distance with cell center $\xx_k$, i.e.
\begin{equation} \label{eq.min_dist}
	\xx_i = \arg \min_{\cell{j} \in \tilde{\mathcal{T}}_{\star}} | \xx_k - \xx_j |.  
\end{equation}   
At line $k$ of linear system \eqref{eq.alg}, the discretization of the differential problem is substituted by
\begin{equation}
	\label{eqn.rec_fringe}
	\Phi_k = \Pol{i} (\Phi(\xx_k)) = \Phi_i + (\tilde{\zz}_i^P(\xx_k))^T \tilde{\boldsymbol{\beta}}_i.
\end{equation}
Consequently, by inserting the definition of the expansion coefficients \eqref{eq.pol_coeff} in \eqref{eqn.rec_fringe}, the $k$-th line of the algebraic problem becomes
\begin{equation} \label{eq.k-th}
	\Phi_i + (\tilde{\zz}_i^P(\xx_k))^\top \tilde{\mathbf{Z}}_i^{P,\dagger} (\boldsymbol{\Phi}_j - \boldsymbol{\Psi}_i) - \Phi_k = 0. 
\end{equation}
Since vector $\hat{\boldsymbol{\Phi}}$ is unknown, the $k$-th line of system \eqref{eq.alg} given by \eqref{eq.k-th} is arranged in order to have $\mathbf{f}_k = 0$, $\mathbf{L}_{kk} = 1$, $\mathbf{L}_{ki} = -1 + (\tilde{\zz}_i^P(\xx_k))^\top \tilde{\mathbf{Z}}_i^{P,\dagger} \mathbf{1}$ and $\mathbf{L}_{kj} = (\tilde{\zz}_i^P(\xx_k))^\top \tilde{\mathbf{Z}}_i^{P,\dagger}$ for any $\cell{j} \in \Ste{i}$, with $\mathbf{1} \in \R^{N_s}$ the identity vector.

\begin{remark}[One-shot extrapolation]
	The proposed trick \eqref{eq.k-th} allows to find the solution of an algebraic problem and, at the same time, to extrapolate values on fringe cells without using an iterative procedure (e.g., Schwartz approach) or trying to find a proper discretization on fringe cells as done in \cite{bergmann2022ader}, for instance. 
\end{remark}  

\begin{remark}[Identity property on aligned mesh blocks]
	If the distance minimizer $\xx_i$ for a fringe cell $\cell{k}$ coincides with the cell center $\xx_k$ (i.e. $| \xx_i - \xx_k | = 0$ in \eqref{eq.min_dist}), we obtain the identity relation $\Phi_i = \Phi_k$ in \eqref{eqn.rec_fringe}. In fact, since the polynomial coefficients are chosen in order to impose exactly $\hat{\beta}_{i0}=\Phi_i$, this implies that if the distance between cells $\cell{i}$ and $\cell{k}$ is zero, the vector of basis functions $\tilde{\zz}_i^P(\xx_k)$ vanishes. In this way, if the mesh blocks are chosen in order to have zero shifting (namely, as there was one mesh block because all cells in the overlapping zone perfectly coincide), the algebraic problem automatically looks for a solution as it was only one block mesh.
\end{remark}

\subsection{High-order semi-implicit IMEX scheme}
To achieve higher order accuracy in time, the class of semi-implicit Implicit-Explicit (IMEX) Runge-Kutta methods is employed \cite{boscarino2016high}. An IMEX Runge-Kutta scheme is a multi-step method characterized by two $s \times s$ triangular matrices: the explicit matrix, referred to as $\tilde{\mathbf{A}} = [\tilde{a}_{ij}]$, where $ \tilde{a}_{ij} = 0$ for $j \geq i$, and the implicit matrix, referred to as $\mathbf{A} = [a]_{ij}$, with $a_{ij} = 0$ for $j > i$ since we use diagonally implicit Runge-Kutta schemes. The number of implicit Runge-Kutta stages is given by $s$. Additionally, the scheme is defined by weight vectors $\tilde{\mathbf{b}}^\top$ and $\mathbf{b}^\top$ in $\R^s$. These matrices and vectors are usually presented in the form of explicit (on the left) and implicit (on the right) Butcher tableaux
\begin{center}
	\begin{tabular}{c | c}
		$\tilde{\mathbf{c}}$ & $\tilde{\mathbf{A}}$ \\ \hline
		& $\tilde{\mathbf{b}}^\top$
	\end{tabular}
	\hspace{1.0cm}
	\begin{tabular}{c | c}
		$\mathbf{c}$ & $\mathbf{A}$ \\ \hline
		& $\mathbf{b}^\top$
	\end{tabular}
\end{center} 
where vectors $\tilde{\mathbf{c}}$ and $\mathbf{c}$ in $\R^s$ are the sum of lines of the explicit and implicit matrices $\tilde{\mathbf{A}}$ and $\mathbf{A}$, respectively. Namely,
\begin{equation} \label{eq.c}
	\tilde{c}_i = \sum_{j = 1}^s \tilde{a}_{ij} \quad \text{and} \quad c_i = \sum_{j = 1}^s a_{ij}, \quad i = 1, \dots, s,
\end{equation}
which provide the time levels of the Runge-Kutta stages. Following \cite{boscarino2016high}, the semi-discrete first order scheme \eqref{eq.ustar}-\eqref{eq.pressure_int} presented in Section \ref{sec.IMEX1} for velocity $\uu^{n+1}$ and pressure field $p^{n+1}$, as well as for the position vector $\xx^{n+1}$ through the motion equation \eqref{eq.space_problem}, can be resumed in the abstract formulation
\begin{equation} \label{eq.abstractQ}
	\frac{\partial \hat{\QQ}}{\partial t} = \HH(t, \EX{\hat{\QQ}}(t), \IM{\hat{\QQ}}(t)),
\end{equation}
where the state vector $\hat{\QQ}=(\QQ,\xx)$ includes both physical quantities $\QQ$ and the evolving frame $\xx$, while the flux $\HH$ accounts for both implicit and explicit terms. The direct dependency of $\HH$ on $t$ remarks that system \eqref{eq.abstractQ} is not autonomous because the space is subject to the deformation due to the mesh velocity $\w^{n+1}$ in \eqref{eq.space_problem}. A partitioned system for $\hat{\QQ} = (\EX{\hat{\QQ}}, \IM{\hat{\QQ}})$ is then defined as follows
\begin{equation} \label{eq.partionedQ}
	\begin{aligned}
		\frac{\partial \EX{\hat{\QQ}}}{\partial t} &= \HH(t, \EX{\hat{\QQ}}(t), \IM{\hat{\QQ}}(t)), \\
		\frac{\partial \IM{\hat{\QQ}}}{\partial t} &= \HH(t, \EX{\hat{\QQ}}(t), \IM{\hat{\QQ}}(t)),
	\end{aligned}
\end{equation}
which is integrated from initial time $t^n$ over the time interval $\dt$ using a partitioned Runge-Kutta method as 
\begin{equation} \label{eq.partitionedQ2}
	\begin{aligned}
		k_i &= \HH \left( t^n + \tilde{c}_i \dt, \QQ^n + \dt \sum_{j = 1}^s \tilde{a}_{ij} k_j, \QQ^n + \dt \sum_{j = 1}^s a_{ij} \ell_j \right) \\
		\ell_i &= \HH \left( t^n + c_i \dt, \QQ^n + \dt \sum_{j = 1}^s \tilde{a}_{ij} k_j, \QQ^n + \dt \sum_{j = 1}^s a_{ij} \ell_j \right)
 	\end{aligned}, \quad 1 \leq i \leq s,
\end{equation}
where the initial implicit and explicit states are given by $\EX{\hat{\QQ}}^n = \IM{\hat{\QQ}}^n = \hat{\QQ}^n$. For a general IMEX Runge-Kutta scheme, both systems in \eqref{eq.partitionedQ2} should be solved providing two sets of fluxes, namely $k_i$ and $\ell_i$. Since the original system \eqref{eq.abstractQ} is not autonomous, the only way for integrating the system with only one set of fluxes of the partitioned system \eqref{eq.partitionedQ2} is to consider Runge-Kutta schemes with Butcher tableaux defined by explicit and implicit matrices whose lines have the same summation, i.e., $\tilde{\mathbf{c}} = \mathbf{c}$, meaning that the stages are defined at the same intermediate time levels, see \cite{boscarino2016high} for further details. For this reason, the numerical approach devised in the sequel only works under this hypothesis. 

At this point, we can consider only one evaluation of \eqref{eq.partitionedQ2} by computing only one set of stage fluxes, that is
\begin{equation} \label{eq.partitionedQ3}
	k_i = \HH \left( t^n + c_i \dt, \QQ^n + \dt \sum_{j = 1}^s \tilde{a}_{ij} k_j, \QQ^n + \dt \sum_{j = 1}^s a_{ij} k_j \right), \quad 1 \leq i \leq s.
\end{equation}
For any $i$-th stage, the implicit problem \eqref{eq.partitionedQ3} is solved through the following scheme
\begin{equation} \label{eq.implicitQ}
	\begin{aligned}
		\EX{\hat{\QQ}}^i &= \hat{\QQ}^n + \dt \sum_{j = 1}^{i-1} \tilde{a}_{ij} k_j , \\
		\IM{\tilde{\hat{\QQ}}}^i &= \hat{\QQ}^n + \dt \sum_{j = 1}^{i-1} a_{ij} k_j , \\
		k_i &= \HH \left( t^n + c_i \dt, \EX{\hat{\QQ}}^i, \IM{\tilde{\hat{\QQ}}}^i + k_i \, a_{ii} \dt \right).
	\end{aligned}
\end{equation}
The implicit problem in \eqref{eq.implicitQ} for $k_i$ is characterized by a time step $a_{ii} \dt$ and it considers a space-time evolution given by $t^n + c_i \dt$ at any implicit iteration stage $i$. The final numerical solution is then assembled by
\begin{equation} \label{eq.update}
	\hat{\QQ}^{n+1} = \hat{\QQ}^n + \dt \sum_{i = 1}^s b_i k_i .
\end{equation}

Now, we provide all the details related to the evaluation of any implicit Runge-Kutta stage of index $i$, with $1 \leq i \leq s$, following the steps of the first order semi-discrete scheme presented at the end of Section \ref{sec.IMEX1}.
\begin{enumerate}
	\item We start by solving implicitly the motion equation \eqref{eq.motion} as
	\begin{equation} \label{eq.spacetime}
		\begin{aligned}
			\EX{\xx}^{(i)} &= \xx^n + \dt \sum_{j = 1}^{i-1} \tilde{a}_{ij} \mathbf{w}_j, \\
			\IM{\tilde{{\xx}}}^i &= \xx^n + \dt \sum_{j = 1}^{i-1} {a}_{ij} \mathbf{w}_j, \\
			\mathbf{w}_i &= \mathbf{w}( t^n + c_i \dt, \tilde{\xx}^{(i)}_I + \mathbf{w}_i a_{ii} \dt ),
		\end{aligned}
	\end{equation}
where the computation of the implicit fluxes $\mathbf{w}_i$ corresponds to the solution of \eqref{eq.space_problem}. 
    The new mesh configuration at time $t=t^n+c_i\dt$ is then given by 
    \begin{equation}
    	\IM{\xx}^{(i)} = \xx^n + \dt \sum_{j = 1}^{i} a_{ij} \mathbf{w}_j.
    \end{equation}
 Moreover, at each stage $i$, we define 
	\begin{equation}
		\old{\xx}^{(i)} = \IM{\xx}^{(i-1)} \quad \text{ with } \quad \old{\xx}^{(0)} = \xx^n,
	\end{equation}
	that provides the implicit spatial configuration relative to the previous Runge-Kutta stage $i-1$. In the sequel, given a generic quantity $\Phi$, the subscripts $\EX{\Phi}$, $\IM{\Phi}$ and $\old{\Phi}$ refer to the quantity $\Phi$ evaluated on the explicit ($\EX{\xx}$), implicit ($\IM{\xx}$) and old ($\old{\xx}$) mesh configuration according to the above definitions.
	
	\item The intermediate velocity field $\uu^{*(i)}$ is obtained by solving \eqref{eq.ustar}. Let us recall the discrete Laplacian, convective and gradient operators given by $\mathbb{K}_h$ in \eqref{eq.laplacian_approx}, $\mathbb{F}_h$ in \eqref{eqn.convOp} and $\mathbb{G}_h$ in \eqref{eq.disc_gradp_divu}, respectively. The momentum equation \eqref{eq.ustar} is discretized for the IMEX stage $i$ as follows:	
	\begin{equation} \label{eqn.ustar_imex}
	  \IM{|\omega|} \, \IM{\uu}^{*(i)} - \frac{a_{ii} \dt}{\Rey} \left[\mathbb{K}_h(\uu^{*(i)})\right]_I = \EX{|\omega|} \, \EX{\tilde{\uu}}^{*(i)} - a_{ii} \dt \left[\mathbb{F}_h(\uu^{(i)},\w^{(i)})\right]_E - a_{ii} \dt \left[\mathbb{G}_h(p^{(i)})\right]_{\text{old}}.
	\end{equation}
Notice that the pressure gradient is evaluated at the old implicit configuration. This is a crucial detail since the pressure is only obtained by the projection-correction equation \eqref{eq.pressure_int} that is indeed solved at the implicit stages. As such, we only have information about the pressure related to the mesh configuration $\old{\xx}$. 

	\item Once the discrete intermediate velocity field $\IM{\uu}^{*(i)}=(u^{*(i)},v^{*(i)})$ is computed through \eqref{eqn.ustar_imex}, we can solve the pressure equation \eqref{eq.pressure_int} on the implicit mesh configuration as
	\begin{equation} \label{eq.alg_pressure}
		a_{ii} \dt \left[\mathbb{K}_h(p^{(i)})\right]_I = \left[\mathbb{G}_h(u^{*(i)}) + \mathbb{G}_h(v^{*(i)}) \right]_I + a_{ii} \dt \left[\mathbb{K}_h(p^{(i)})\right]_{\text{old}}.
	\end{equation}

	\item Finally, the divergence-free velocity is updated with the previously computed implicit pressure $p^{(i)}_I$ as
	\begin{equation} \label{eq.alg_update}
		\IM{\uu^{(i)}} = \IM{\uu^{*(i)}} - a_{ii} \dt \left[\mathbb{G}_h(p^{(i)})\right]_I.
	\end{equation} 
\end{enumerate}

\begin{remark}[Evolving space-time mesh configurations]
	In both \eqref{eqn.ustar_imex} and \eqref{eq.alg_pressure}, the terms related to the pressure refer to the old implicit space-time mesh configuration. This is due to the fact the IMEX fractional-step method provides the pressure at each stage by solving the projection-correction equation \eqref{eq.alg_pressure}, but no evolution equation is available for the pressure. Thus, the information on the pressure always refers to the old implicit configuration but not to the current explicit space-time configuration needed in \eqref{eq.implicitQ} for the computation of the intermediate states. For the velocity it is different, since we have an evolution equation for the velocity eventually given by \eqref{eq.alg_update}. As a matter of fact, in general it holds that $\old{\xx}^{(i)} = \IM{\xx}^{(i-1)} \neq \EX{\xx}^{(i)}$. On the other hand, the contrary is true for first order Runge-Kutta schemes and for not evolving frame configurations, namely, when all space-time configurations coincide with the space-time configuration at $t = 0$, hence no mesh motion is considered.  
\end{remark}

\begin{remark}[Butcher tableaux of the semi-implicit IMEX schemes]
An IMEX scheme is described with a triplet $(s,\tilde{s},p)$ which characterizes the number $s$ of stages of the implicit method, the number $\tilde{s}$ of stages of the explicit method and the order $p$ of the resulting scheme. The employed Runge-Kutta schemes are the \textit{forward-backward Euler (1,1,1)} and the \textit{L-stable two-stage Diagonally Implicit ARS (2,2,2)}, for achieving first and second order of accuracy, respectively. The reader is addressed to \cite{ascher1997implicit} for their derivation. The Butcher tableaux for these two IMEX Runge-Kutta schemes read:
\begin{itemize}
	\item Euler (1,1,1)
	\begin{equation}
		\begin{array}{c|c}
			0 & 0 \\ \hline & 1
		\end{array} \qquad
		\begin{array}{c|c}
			1 & 1 \\ \hline & 1
		\end{array}
		\label{eq.imex_euler}
	\end{equation}

\item ARS (2,2,2)
\begin{equation}
	\begin{array}{c|cc}
		\gamma & \gamma & 0 \\ 1 & \delta & 1 - \delta \\ \hline & \delta & 1 - \delta
	\end{array} \qquad
	\begin{array}{c|cc}
		\gamma & \gamma & 0 \\ 1 & 1 - \gamma & \gamma \\ \hline & 1-\gamma & \gamma
	\end{array}
	\label{eq.imex_ars}
\end{equation}
The parameters in \eqref{eq.imex_ars} are $\gamma = 1 - \sqrt{2}/2$ and $\delta = 1 - 1/(2 \gamma)$.
\end{itemize}   

\end{remark}

\subsection{Evolving overset configurations}
During the simulation, as the foreground mesh moves, the background mesh undergoes changes in both the overlapping region and the hole. Let us consider a background cell, denoted as $\cell{i}$, located in the neighborhood of the overlapping area. Between two consecutive times or stages of the implicit Runge-Kutta time stepping, one of the following three scenarios may occur:
\begin{enumerate}
	\item The cell $\cell{i}$ is present at the current time but it disappears at the next time level because the hole fully covers it.
	\item The cell $\cell{i}$ is not present at the current time, but it appears at the next time level as the hole moves away from it.
	\item The configuration of the overlapping zone remains unchanged with respect to the cell $\cell{i}$. Consequently, the cell is present at both current and next times.
\end{enumerate}

In the first case, the algorithm processes the vanishing cell to compute the fluxes required by updating the neighboring cells, and at the next time, the cell $\cell{i}$ with its data is simply removed. 

For the second scenario, the information from the current time is missing, and it is necessary to extrapolate the same information in order to advance the solution of the new cell $\cell{i}$ to the next time level. This process is performed using the same approach as for a fringe cell. Specifically, if the new cell is born, it is certainly a background fringe cell. So we search for a foreground cell $\cell{j}$ that minimizes the Euclidean distance between the cell centers $\xx_i$ and $\xx_j$ according to \eqref{eq.min_dist} with $\tilde{\mathcal{T}}_{\star}=\mathcal{T}_{\text{bg}}$. Next, an interpolating polynomial $\Pol{j}(\Phi)$ (for $\Phi$ any relevant physical variable) centered on $\xx_j$ is constructed, and the data that is stored at the cell center of $\cell{i}$ is the evaluation of that polynomial at $\xx_i$. This preserves the second order of spatial accuracy of the scheme.

The third case is straightforward. 

\section{Numerical results} \label{sec.numerical_results}

The novel numerical schemes are applied to several test cases in order to assess convergence, stability and accuracy properties. Whenever possible, the different benchmarks are compared against exact or numerical solutions. The quantitative analysis is carried out in $L^2$-norm. In particular, let $\Phi : \Omega \rightarrow \R$ be a specific function defined in the computational domain $\Omega$, its $L^2$-norm reads
\begin{equation} \label{eq.L2}
	\| \Phi \|_{L^2} = \sqrt{ \int_\Omega \Phi^2 \, \diff \xx }.
\end{equation}

For all tests, the CFL number in \eqref{eq.CFL} is assumed equal to 0.9. It is important to emphasize that the time step remains independent of the fast scales in the problem being addressed, as these scales are discretized implicitly (pressure and viscosity terms). Additionally, when the velocity field is initialized with zero, the determination of the first time step follows the CFL condition commonly used in fully explicit schemes, thereby incorporating the eigenvalues of the entire system of governing equations \cite{boscheri2023new}.

This section is organized as follows. Firstly, we test the convergence of the second order finite volume IMEX Runge-Kutta scheme. We also demonstrate that the schemes fulfill the free-stream preservation property, also known as Geometric Conservation Law up to machine accuracy on moving Chimera meshes. Successively, we numerically prove that the scheme is precise at zero-machine if the exact solution of the problem is at most a polynomial of degree two through the analytical solution of the Poiseuille flow. Different overset configurations are employed in the lid-driven cavity test for ensuring that no loss of properties are given by the overlapping zone. Finally, different benchmarks on cylinders are proposed for analyzing the accuracy of the method.

If it is not specified, the method is intended to be second-order convergent, thus we use the polynomial reconstruction $\mathcal{R}^{P,Q}$ in space and the ARS (2,2,2) scheme in time. For the foreground mesh, its original location is specified (i.e. at time $t = 0$). Successively, the overlapping zone is built by rimming the outermost foreground cell layer with 5 layers of cells similar in size to those of the background cell.

\subsection{Converge rate on Taylor-Green vortexes}
The Taylor-Green test describes an unsteady inviscid flow of a vortex in the domain $\Omega = [-\pi, \pi]^2$ with periodic boundary conditions for both velocity field and pressure. The exact solution for this problem is
\begin{equation} \label{eq.TGV}
	p = -\frac{e^{-4 t / \Rey}}{4} \left( \cos( 2x ) + \cos( 2y ) \right) , \qquad 
	\uu = \begin{bmatrix}
		u \\ v
	\end{bmatrix} = \begin{bmatrix}
      \phantom{-}\sin( x ) \cos(y)  \\
      - \cos( x ) \sin( y) 
\end{bmatrix} e^{-2 t / \Rey}.
\end{equation}  
The initial condition is defined by the solution \eqref{eq.TGV} evaluated at time $t = 0$.

A convergence analysis is conducted for the second order IMEX Runge-Kutta scheme with second order polynomial reconstruction. The analysis is carried out on four refined computational meshes, each of them with three different Reynolds numbers $\Rey = \{ 10^1, 10^3, 10^6 \}$. Each mesh is identified by its characteristic length $\bar{h}$, which is determined as the maximum element size in the computational grid. The errors are measured in $L^2$-norm for the $x$-direction velocity $u$ and for the pressure $p$ at the final time $t_f = 0.2$. The foreground mesh is originally collated in the subdomain $[-1, 1]^2$ as a Cartesian mesh with the same cell size of the background mesh. It is prescribed a space-time dependent velocity 
\begin{equation}
	\w = [w_x,w_y]^\top = 0.2 e^{t - t_f} [ \sin(x) \cos(y), \cos(x) \sin(y) ]^\top.
\end{equation}

\begin{figure}[!htbp]
	\begin{center}
		\begin{tabular}{cc}
			\includegraphics[width=0.49\textwidth]{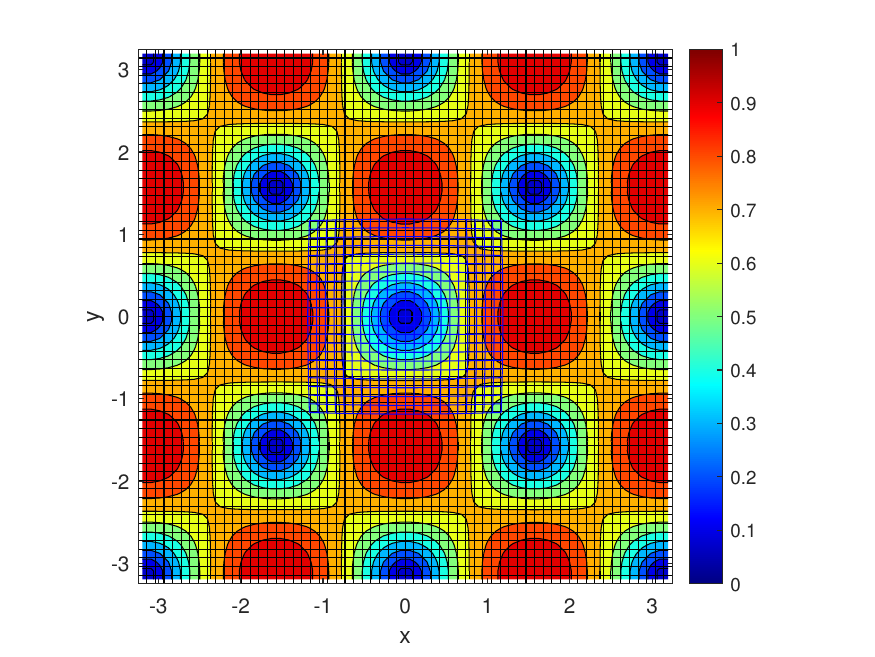} &
			\includegraphics[width=0.49\textwidth]{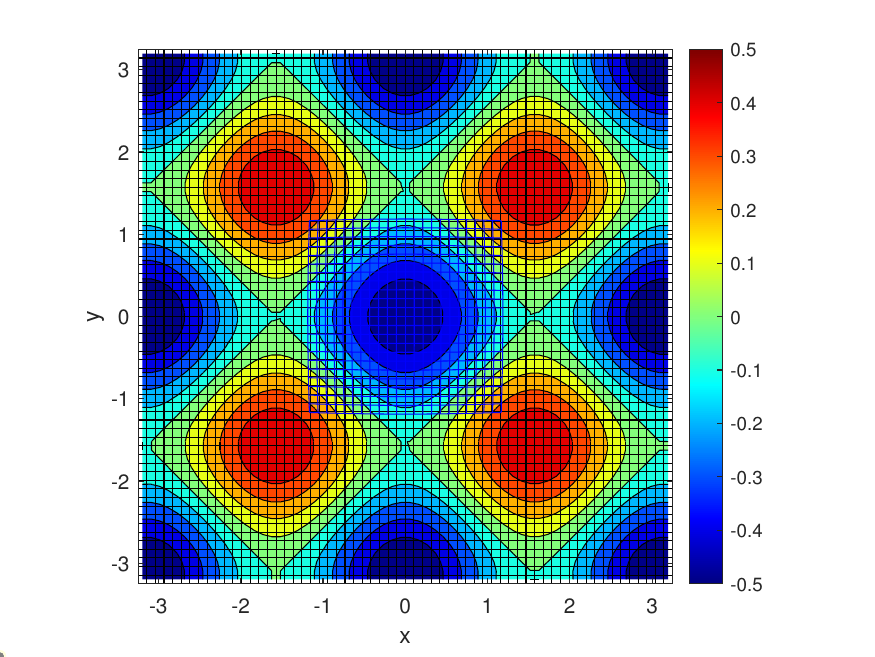}
		\end{tabular}
	\end{center}
	\caption{Taylor-Green vortex at final time $t_f = 0.2$ at Reynolds number $\Rey = 10^3$. Left: velocity magnitude $| \uu |$. Right: pressure field $p$. The overset configuration is defined by a black background grid and a blue foreground mesh. }
	\label{fig.TGV_magn_pres}
\end{figure}

\begin{figure}[!htbp]
	\begin{center}
		\begin{tabular}{cc}
			\resizebox{.49\linewidth}{!}{\begin{tikzpicture}
				\begin{axis}[
					xlabel = {$x$},
					ylabel = {$v,p$},
					xmin = -3.14, xmax = 3.14,
					ymin = -1, ymax = 1,
					line/.style={thick},
					legend style = {at={(axis cs:3.1,1.1)}, anchor=north east},
					grid=both,
					grid style={line width=.1pt, draw=gray!10}]
					
					\addplot[blue, mark=o, only marks] table {TGV_vbg1.dat}; \addlegendentry{$v$ (bg)}
					\addplot[red, mark=x, only marks] table {TGV_vfg.dat}; \addlegendentry{$v$ (fg)}
					\addplot[black, solid] table {TGV_vbg1_ex.dat}; \addlegendentry{$v$ (exact)}
					
					\addplot[blue, mark=+, only marks] table {TGV_pxbg1.dat}; \addlegendentry{$p$ (bg)}
					\addplot[blue, mark=*, only marks] table {TGV_pxfg.dat}; \addlegendentry{$p$ (fg)}
					\addplot[black, dashed] table {TGV_pxbg1_ex.dat}; \addlegendentry{$p$ (exact)}
					
					\addplot[blue, mark=o, only marks] table {TGV_vbg2.dat}; 
					\addplot[black, solid] table {TGV_vbg2_ex.dat}; 
					\addplot[black, solid] table {TGV_vfg_ex.dat}; 
					
					\addplot[blue, mark=+, only marks] table {TGV_pxbg2.dat}; 
					\addplot[black, dashed] table {TGV_pxbg2_ex.dat};
					\addplot[black, dashed] table {TGV_pxfg_ex.dat};
					
				\end{axis}
			\end{tikzpicture} } &
			\resizebox{.49\linewidth}{!}{ \begin{tikzpicture}
				\begin{axis}[
					xlabel = {$y$},
					ylabel = {$u,p$},
					xmin = -3.14, xmax = 3.14,
					ymin = -1, ymax = 1,
					line/.style={thick},
					legend style = {at={(axis cs:-3.1,1.1)}, anchor=north west},
					grid=both,
					grid style={line width=.1pt, draw=gray!10}]
					
					\addplot[blue, mark=o, only marks] table {TGV_ubg1.dat}; \addlegendentry{$u$ (bg)}
					\addplot[red, mark=x, only marks] table {TGV_ufg.dat}; \addlegendentry{$u$ (fg)}
					\addplot[black, solid] table {TGV_ubg1_ex.dat}; \addlegendentry{$u$ (exact)}
					
					\addplot[blue, mark=+, only marks] table {TGV_pybg1.dat}; \addlegendentry{$p$ (bg)}
					\addplot[blue, mark=*, only marks] table {TGV_pyfg.dat}; \addlegendentry{$p$ (fg)}
					\addplot[black, dashed] table {TGV_pybg1_ex.dat}; \addlegendentry{$p$ (exact)}
					
					\addplot[blue, mark=o, only marks] table {TGV_ubg2.dat}; 
					\addplot[black, solid] table {TGV_ubg2_ex.dat}; 
					\addplot[black, solid] table {TGV_ufg_ex.dat}; 
					
					\addplot[blue, mark=+, only marks] table {TGV_pybg2.dat}; 
					\addplot[black, dashed] table {TGV_pybg2_ex.dat};
					\addplot[black, dashed] table {TGV_pyfg_ex.dat};
					
				\end{axis}
			\end{tikzpicture} }
		\end{tabular}
	\end{center}
	\caption{Taylor-Green vortex at final time $t_f = 0.2$ with Reynolds number $\Rey = 10^3$. Left: 1D cuts along the $x$-axis of the $y$-component $v$ of the velocity field $\uu$ and of the pressure $p$ against the exact solution. Right: 1D cuts along the $y$-axis of the $x$-component $u$ of the velocity field $\uu$ and of the pressure $p$ against the exact solution. The two background zones for $x,y \in [-\pi, 1.1] \cup [1.1, \pi]$ consist of 30 equidistant points per zone. The foreground zone is defined by 20 equidistant points in interval $[-1.1, 1.1]$.}
	\label{fig.TGV_cuts}
\end{figure}
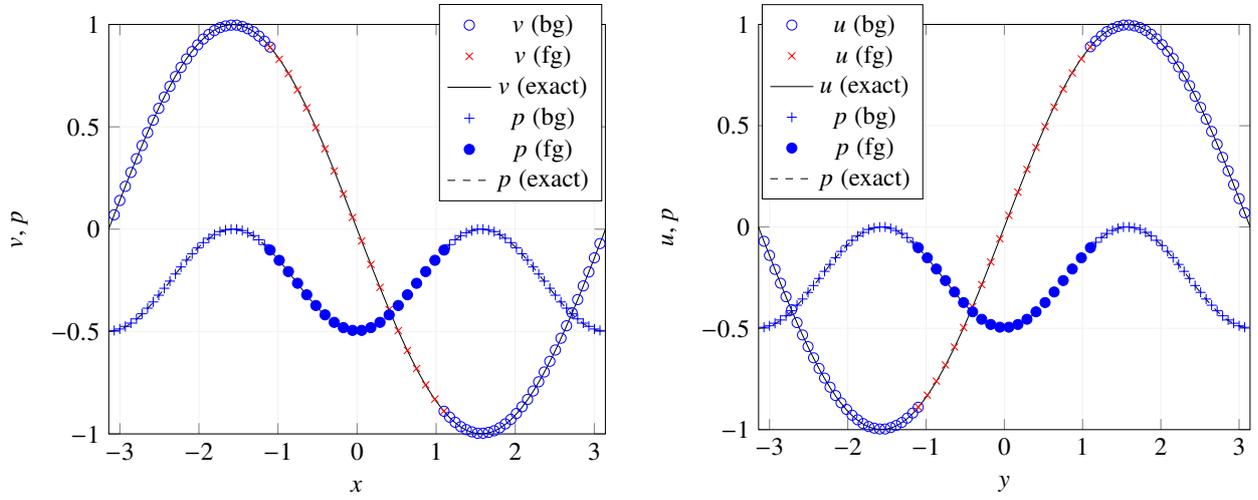

Table \ref{tab.TGC_conv} presents a summary of the results and reports the achieved convergence rates. Remarkably, the method demonstrates the correct order of accuracy while maintaining both asymptotic preservation and accuracy, with the order of accuracy remaining consistent across different viscosity values. We underline that the time step is the same independently of the viscosity coefficient, hence asymptotic accuracy is numerically observed. Figure \ref{fig.TGV_magn_pres} depicts the final Chimera configuration with the velocity magnitude $|\uu|$ (left) and pressure field $p$ (right). Along the axis $x = 0$ and $y = 0$, we report the comparison of the achieved numerical solution with the exact solution \eqref{eq.TGV} in Figure \ref{fig.TGV_cuts} for $\Rey=10^3$.

\begin{table}[!htbp]
	\label{tab.TGC_conv}
	\caption{Convergence rate studies for the Taylor-Green vortex. The study is conducted for three different Reynolds number $\Rey \in \{10^1, 10^3, 10^6\}$ with a fully second order reconstruction in space and a second order IMEX Runge-Kutta scheme in time. The errors are measured in $L^2$-norm and refer to the $x$-component $u$ of velocity field $\uu$ and pressure $p$ with respect to the maximum cell size $\bar{h}$ at final time $t_f = 0.2$.}
	\centering
	\begin{tabular}{llllll}
		\toprule
		$\Rey$                  & $\bar{h}$ & $\| u \|_{L^2}$ & $\mathcal{O}(u)$ & $ \| p \|_{L^2} $ & $\mathcal{O}(p)$ \\
		\midrule
		\multirow{4}{*}{$10^1$} & 3.4261E-1 & 1.2614E-2       & $-$              & 8.2517E-2         & $-$              \\
		& 1.7181E-1 & 4.1789E-3       & 1.60             & 1.7119E-2         & 2.27             \\
		& 1.1460E-1 & 1.6332E-3       & 2.32             & 7.3411E-3         & 2.10             \\
		& 8.5969E-2 & 8.5412E-4       & 2.26             & 3.5839E-3         & 2.49             \\
		\midrule
		\multirow{4}{*}{$10^3$} & 3.4261E-1 & 1.4294E-2       & $-$              & 1.2277E-1         & $-$              \\
		& 1.7181E-1 & 4.9396E-3       & 1.53             & 1.3884E-2         & 3.15             \\
		& 1.1460E-1 & 1.9675E-3       & 2.27             & 5.9900E-3         & 2.08             \\
		& 8.5969E-2 & 1.0917E-3       & 2.05             & 3.3840E-3         & 1.99             \\
		\midrule
		\multirow{4}{*}{$10^6$} & 3.4261E-1 & 1.2054E-2       & $-$              & 8.2893E-2         & $-$              \\
		& 1.7181E-1 & 3.0570E-3       & 1.99             & 1.7761E-2         & 2.23             \\
		& 1.1460E-1 & 1.2761E-3       & 2.16             & 7.2470E-3         & 2.21             \\
		& 8.5969E-2 & 7.1376E-4       & 2.02             & 3.9284E-3         & 2.13            \\ 
		\bottomrule
	\end{tabular}
\end{table}

\subsection{Free-stream preservation}
\begin{figure}[!htbp]
	\centering
	\includegraphics[width=0.59\textwidth]{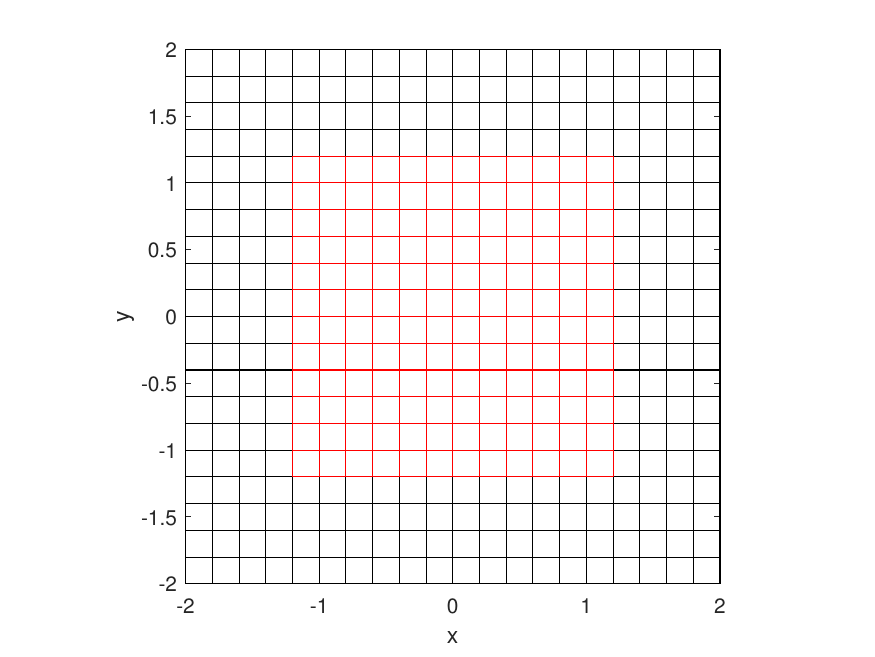}
	\caption{Initial overset configuration for the free-stream preservation test. The foreground mesh (in red) is initially aligned to the background mesh (in black).}
	\label{fig.FSP_initial_configuration}
\end{figure}

Since the proposed method is second order convergent, we expect that the assessed numerical solutions are exact at zero-machine if their exact representation is a polynomial of degree less or equal to two, independently of the evolution of the Chimera mesh. Therefore, we start by showing the performance of the method for the free-stream preservation test, that is also known as Geometric Conservation Law (GCL) in Lagrangian schemes \cite{Maire2007}. We design this test by considering a steady flow of solution $(\uu, p) = (\mathbf{0}, 1)$ for any $t \in [0, 1]$ in the domain $\Omega = [-2, 2]^2$. We employ a coarse background mesh of cell size $\bar{h}=1/5$ and a foreground mesh originally occupying the zone $[-1,1]$ (see Figure \ref{fig.FSP_initial_configuration}). We consider two overset configurations defined by velocities
\[
	\w_1 = 0.3 \cdot [ 1, 1]^\top \quad \text{and} \quad \w_2 = 0.5 \cdot [ y, -x ]^\top,
\]  
accounting for a rigid translation and rotation with respect to the origin of the axis, respectively. The used Reynolds number is $\Rey = 200$.

\begin{figure}[!htbp]
	\begin{center}
		\begin{tabular}{cc}
			\includegraphics[width=0.49\linewidth]{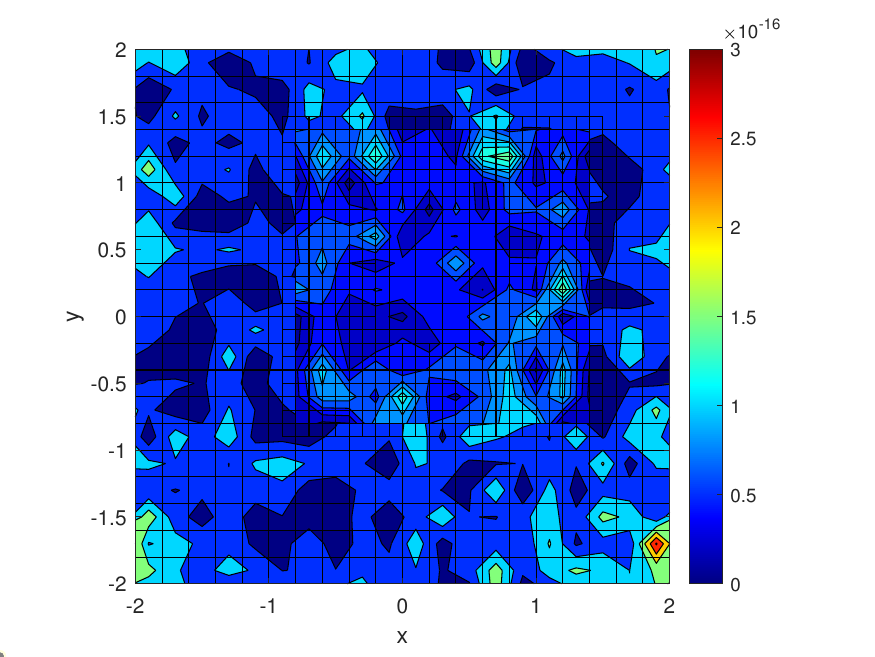} &
			\includegraphics[width=0.49\linewidth]{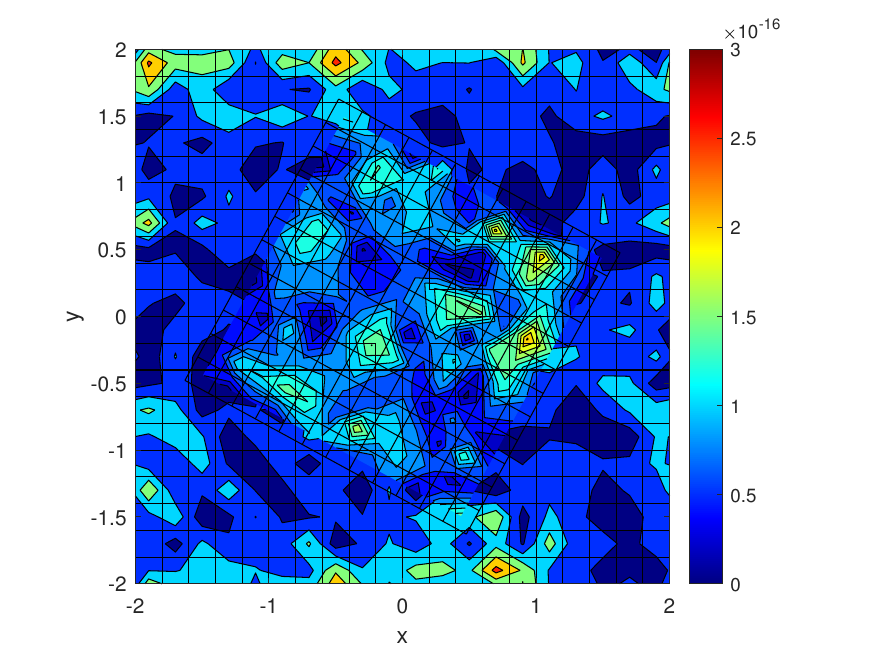}
		\end{tabular}
	\end{center}
	\caption{Velocity magnitude $| \uu |$ for the free-stream preservation test at the final time $t_f = 1$ for foreground mesh translation (left) and rotation (right).}
	\label{fig.FSP_u}
\end{figure}

For both cases, the numerical solution is maintained up to machine accuracy for the whole time interval. In Figure \ref{fig.FSP_u} the pointwise velocity magnitude $| \uu |$ is showed at the final configuration for the translation and rotation cases. The maximum value does not exceed the value of $3 \cdot 10^{-16}$, confirming that the GCL is respected. 

\subsection{Poiseuille flow}
\begin{figure}[!htbp]
	\centering
	\includegraphics[trim={0 3cm 0 3cm},clip, width=0.59\textwidth]{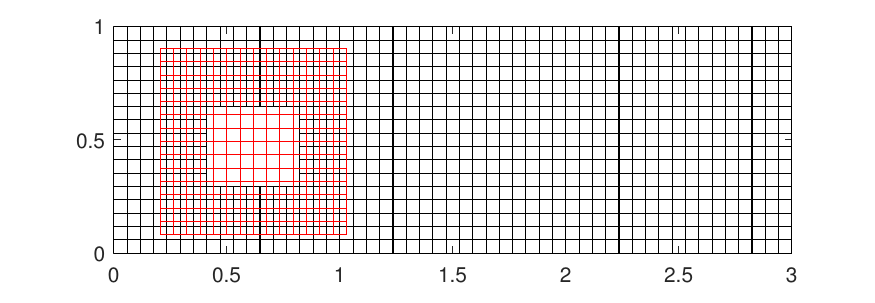}
	\caption{Initial overset configuration for the Poiseuille flow test. The foreground mesh (in red) is initially not aligned to the background mesh (in black).}
	\label{fig.pois_initial_configuration}
\end{figure}

In the previous subsection, we numerically proved the scheme to be precise at zero-machine for a constant solution. Here, we perform a similar test for the Poiseuille flow in the channel $\Omega = [0,3] \times [0,1]$. The exact solution for this test is provided by a steady quadratic velocity profile and linear pressure field, i.e.,
\begin{equation} \label{eq.pois}
	\uu = \begin{bmatrix}
		-y ( y -1 ) \\
		0
	\end{bmatrix}, \quad p(x,y) = - \frac{2x}{\Rey}.
\end{equation}
Also in this case, the Reynolds number is set equal to 200. The initial conditions are defined by \eqref{eq.pois}. The employed background mesh has a cell size of $\bar{h} = 3/50$. The foreground mesh is Cartesian with the same cell size of the background mesh. Initially, it occupies the subdomain $[0.5, 0.75] \times [0.375, 0.625]$ without the overlapping layers (as depicted in Figure \ref{fig.pois_initial_configuration}). The foreground mesh rigidly translates along the $x$-direction with a time-dependent velocity 
\begin{equation}
	\w = -t(t-1) \cdot [1, 0]^2.
\end{equation}

\begin{figure}[!htbp]
	\begin{center}
		\begin{tabular}{cc}
			\includegraphics[trim={0 3cm 0 3cm},clip,width=0.59\linewidth]{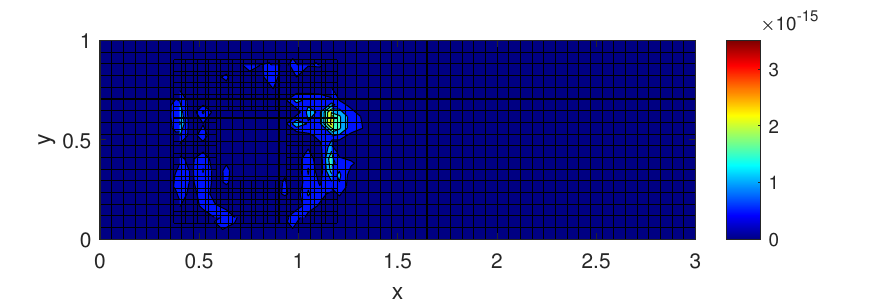} \\
			\includegraphics[trim={0 3cm 0 3cm},clip,width=0.59\linewidth]{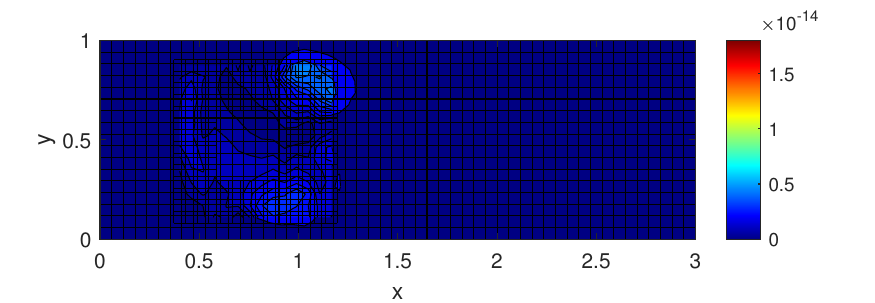}
		\end{tabular}
	\end{center}
	\caption{Pointwise error on $x$-component $u$ of velocity field $\uu$ (top) and pressure $p$ (bottom) for the Poiseuille flow test at final time $t_f = 1$.}
	\label{fig.pois_err}
\end{figure}

In Figure \ref{fig.pois_err} we report the pointwise mismatch between the numerical and exact solution for the $x$-direction velocity $u$ and for the pressure field $p$ at final time $t_f = 1$. Also in this case the second order convergent scheme is precise at zero-machine. As a matter of fact, both velocity and pressure measure an error around $3\cdot 10^{-15}$. 

\subsection{Lid-driven cavity}
\begin{figure}[!htbp]
	\begin{center}
		\begin{tabular}{cc}
			\includegraphics[width=.49\linewidth]{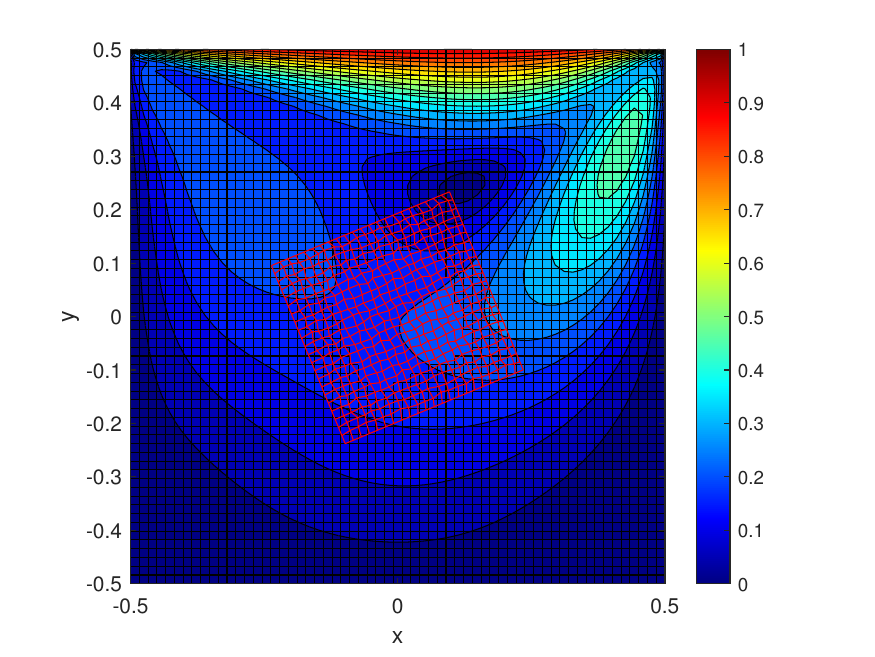} &
			\includegraphics[width=.49\linewidth]{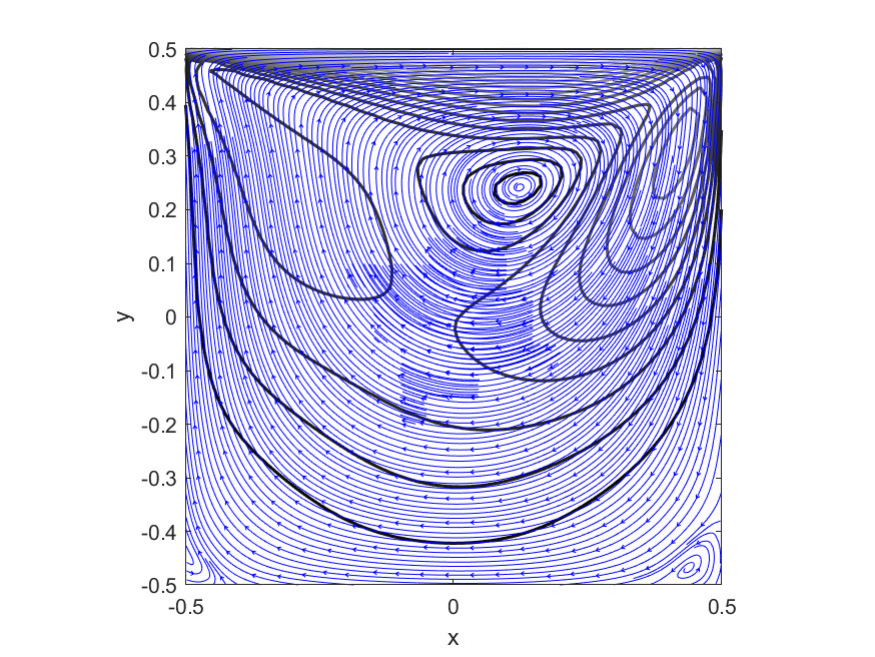} \\
			\includegraphics[width=.49\linewidth]{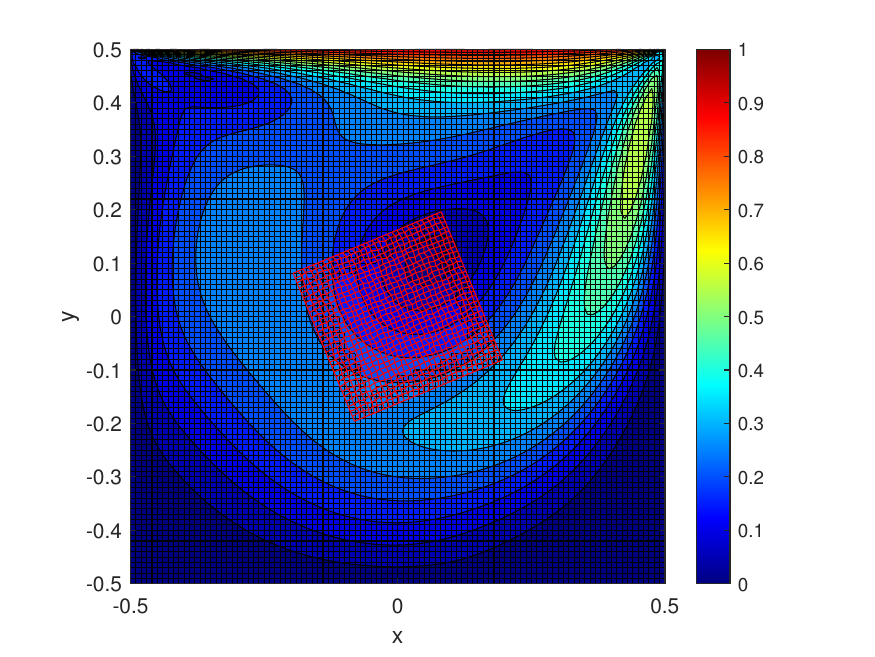} &
			\includegraphics[width=.49\linewidth]{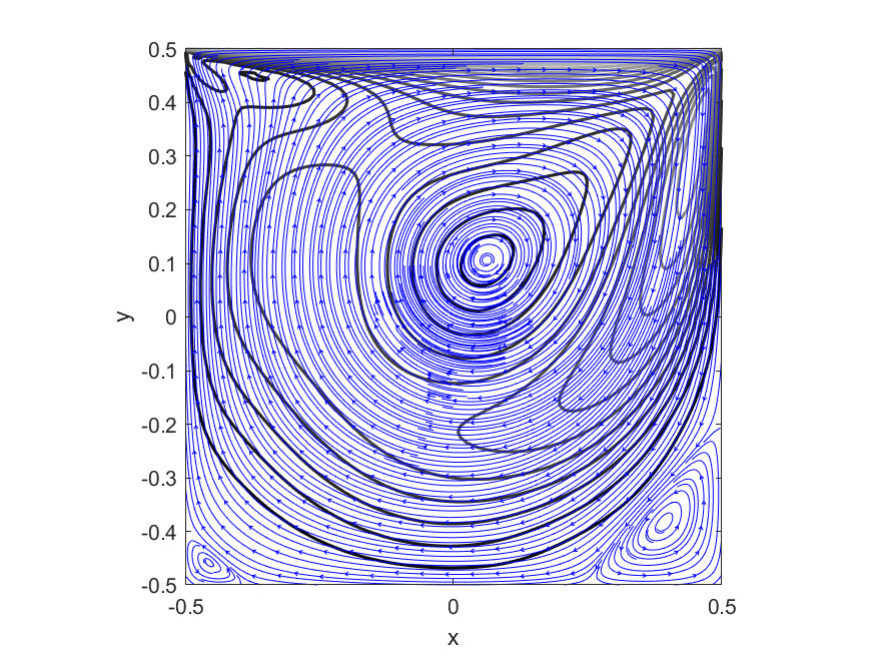}
		\end{tabular}
	\end{center}
	\caption{Velocity magnitude $| \uu |$ for the lid driven cavity test at the final time $t_f = 25$ for Reynolds numbers $\Rey = 100$ (first row) and $\Rey = 400$ (second row). On the left it is possible to see the overset configuration for both cases. The Chimera mesh is steady and an irregular foreground mesh is employed. The mesh at convergence has an average cell size $\bar{h} = 1/60$ and $\bar{h} = 1/100$ for $\Rey = 100$ and $\Rey = 400$, respectively. On the right, the streamlines are superposed to the vorticity lines for both cases.}
	\label{fig.LDC_magn}
\end{figure}

Next, we simulate the classical benchmark of the lid-driven cavity flow within the computational domain $\Omega = [-0.5, 0.5]^2$. The objective in this scenario is to determine the steady hydrodynamics state for a fluid initially at rest, where both pressure and velocity are set to zero at initial time. The domain is bounded by walls on the vertical sides ($x = \pm 0.5$) and along the bottom ($y = -0.5$), while the top side ($y = +0.5$) has an imposed velocity field $\uu = [1, 0]^\top$. We consider two different Reynolds number, namely $\Rey \in \{ 100, 400 \}$. For this test case, the foreground mesh is not moving. It is built by rotating the square defined in the subdomain $[-0.1, 0.1]^2$ by an angle of $\pi/8$. All foreground cells have an average size $\bar{h}$ comparable to the one of the background grid even though their vertexes are randomly displaced of length $0.4 \bar{h}$ with respect to the original Cartesian tesselation (see Figure \ref{fig.LDC_magn}). For $\Rey = 100$, the mesh at convergence is defined by $\bar{h} = 1/60$ and by a total number of active cells $N_c = 4260$. For $\Rey = 400$, the employed mesh at convergence is characterized by $\bar{h} = 1/100$ and $N_c = 10823$. A first-order IMEX Runge-Kutta scheme is adopted since we aim at capturing the steady solution at convergence. The final time is chosen to be $t_f=25$.

\begin{figure}[!htbp]
	\begin{center}
		\begin{tabular}{cc}
			\resizebox{.49\linewidth}{!}{ \begin{tikzpicture}
				\begin{axis}
					[
					xlabel = {$x,y$},
					ylabel = {$v,u$},
					xmin = -0.5, xmax = 0.5,
					ymin = -0.4, ymax = 1.2,
					xtick= {-0.5, -.25, ..., 0.5},
					ytick= {-0.4, -.2, ..., 1.2},
					line/.style={thick},
					legend style = {at={(axis cs:-.45,1.25)}, anchor=north west},
					grid=both,
					grid style={line width=.1pt, draw=gray!10}]
					
					\addplot[blue, solid] table {vbg1_Re100.dat}; \addlegendentry{$v$ (bg)}
					\addplot[blue, dotted] table {vfg_Re100.dat}; \addlegendentry{$v$ (fg)}
					\addplot[black, mark=square, only marks] table {vGGs_100.dat}; ; \addlegendentry{$v$ (reference \cite{ghia1982high})}
					
					\addplot[red, solid] table {ubg1_Re100.dat}; \addlegendentry{$u$ (bg)}
					\addplot[red, dotted] table {ufg_Re100.dat}; \addlegendentry{$u$ (fg)}
					\addplot[black, mark=o, only marks] table {uGGs_100.dat}; ; \addlegendentry{$u$ (reference \cite{ghia1982high})}
					
					\addplot[blue, solid] table {vbg2_Re100.dat}; 
					\addplot[red, solid] table {ubg2_Re100.dat}; 
				\end{axis}
			\end{tikzpicture} } &
			\resizebox{.49\linewidth}{!}{ \begin{tikzpicture}
				\begin{axis}
					[
					xlabel = {$x,y$},
					ylabel = {$v,u$},
					xmin = -0.5, xmax = 0.5,
					ymin = -0.6, ymax = 1.2,
					xtick= {-0.5, -.25, ..., 0.5},
					ytick= {-0.6, -.4, -.2, 0, ..., 1.2},
					line/.style={thick},
					legend style = {at={(axis cs:-.45,1.25)}, anchor=north west},
					grid=both,
					grid style={line width=.1pt, draw=gray!10}]
					
					\addplot[blue, solid] table {vbg1_Re400.dat}; \addlegendentry{$v$ (bg)}
					\addplot[blue, dotted] table {vfg_Re400.dat}; \addlegendentry{$v$ (fg)}
					\addplot[black, mark=square, only marks] table {vGGs_400.dat}; ; \addlegendentry{$v$ (reference \cite{ghia1982high})}
					
					\addplot[red, solid] table {ubg1_Re400.dat}; \addlegendentry{$u$ (bg)}
					\addplot[red, dotted] table {ufg_Re400.dat}; \addlegendentry{$u$ (fg)}
					\addplot[black, mark=o, only marks] table {uGGs_400.dat}; ; \addlegendentry{$u$ (reference \cite{ghia1982high})}
					
					\addplot[blue, solid] table {vbg2_Re400.dat}; 
					\addplot[red, solid] table {ubg2_Re400.dat}; 
				\end{axis}
			\end{tikzpicture} }
		\end{tabular}
	\end{center}
	\caption{Lid-driven cavity flow. Velocity profiles along $x = 0$ and $y = 0$ against reference solutions from \cite{ghia1982high} at the final time $t_f = 25$ for Reynolds numbers $\Rey = 100$ (left) and $\Rey = 400$ (right).}
	\label{fig.LDC_cuts}
\end{figure}
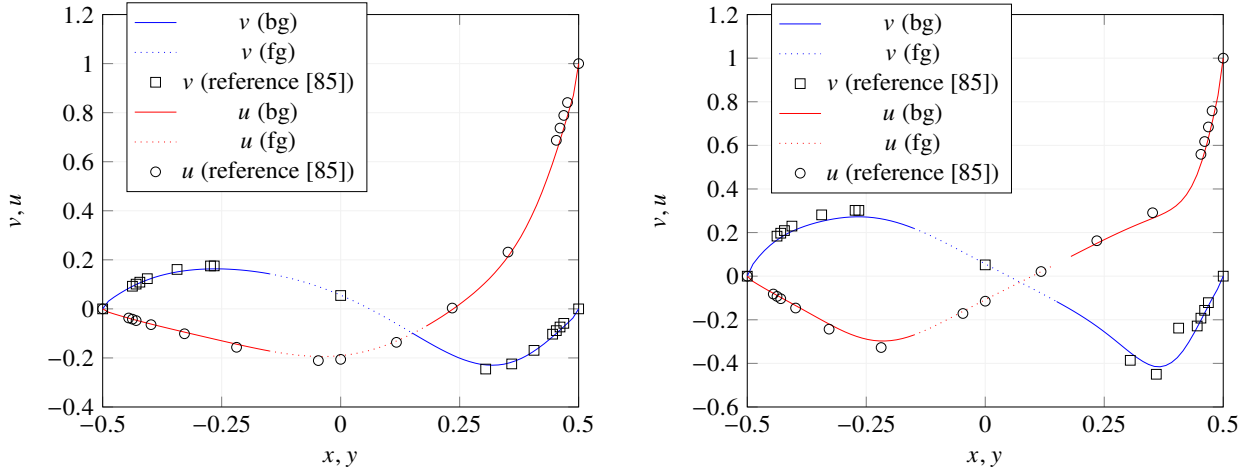

Figure \ref{fig.LDC_magn} depicts the velocity magnitude $|\uu|$, the vorticity and the streamlines at the final time for both considered Reynolds numbers. The streamlines reveal the generation of small vortical flows with an opposite orientation to the main vortex created by the cavity. Additionally, in Figure \ref{fig.LDC_cuts} we present the numerical velocity components $u$ and $v$ along the cuts $x = 0$ and $y = 0$, respectively, and we compare them with the findings reported in \cite{ghia1982high} through direct numerical simulations of incompressible viscous flows. The plots demonstrate a notable agreement between our numerical results and the data available in the literature for both Reynolds numbers.

\subsection{Drag coefficient for steady and moving cylinders} \label{sec.cyl}
In this section, we validate the method by investigating the flow around a cylinder, which may either remain stationary or move. For all tests, the computational domain $\Omega$ is defined by a rectangular channel with a circle $\mathcal{C}$ of center $\xx_c$ and diameter $D$ removed from the circle. In particular, we impose the external left and right boundaries to be inlet and outlet, respectively. Let $\uu_\infty$ be the fluid velocity at the inlet boundary and $\uu_B$ be the body velocity, i.e., the velocity of displacement of the cylinder. We denote the dimensionless stress tensor as $\Tstress(\uu, p)$, which is defined as
\begin{equation} \label{eq.Tstress}
	\Tstress(\uu, p) = -p \I + \frac{1}{\Rey} \left( \nabla \uu + \nabla \uu^T \right).
\end{equation}
The fluid dynamics force $\F_f$ is the integral of the stress tensor \eqref{eq.Tstress} applied to the normal unit vector $\nn_{\mathcal{C}}$ on the surface $\partial \mathcal{C}$:
\begin{equation} \label{eq.Ff}
	\F_f = \oint_{\partial \mathcal{C}} \Tstress(\uu, p) \, \nn_{\mathcal{C}} \, \diff \Gamma_{\mathcal{C}}.
\end{equation}
Force \eqref{eq.Ff} represents the force exerted by the fluid on the cylinder. The aerodynamics drag $C_D$ and lift $C_L$ coefficients are the $x$- and $y$-components of vector $\mathbf{C} = 2 \F_f / (\rho | \uu_c | D)$, respectively. For all test cases, the density $\rho$ of the fluid is set equal to 1 and assumed to be constant. The characteristic velocity $\uu_c$ is equal to $\uu_\infty$ when the cylinder is steady and a flow impacts on it with velocity $\uu_\infty$ at the inlet, otherwise it is equal to the cylinder velocity $\uu_B$ (and the inlet velocity vanishes).

\begin{figure}[!htbp]
	\begin{center}
		\begin{tabular}{cc}
			\includegraphics[width=.49\linewidth]{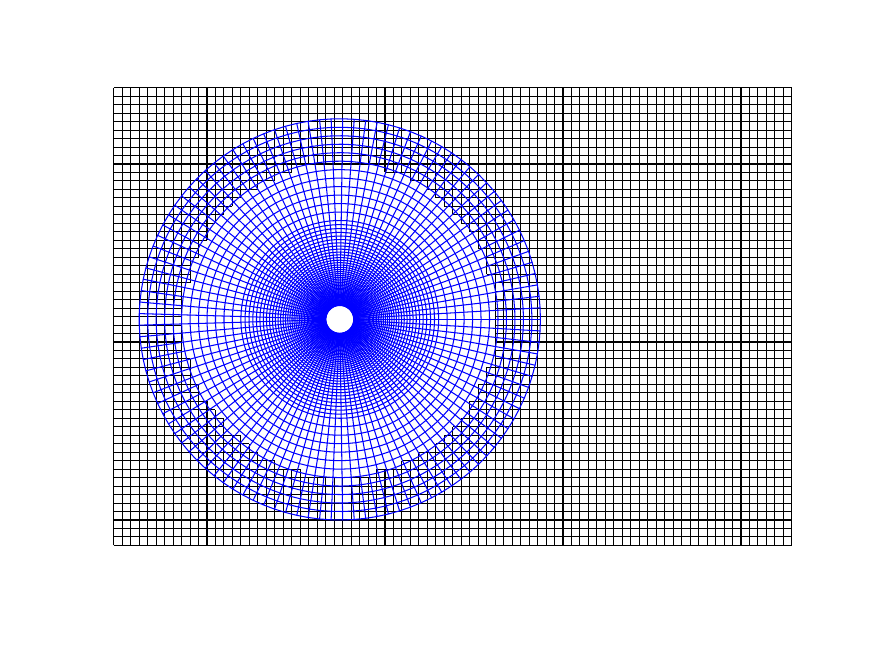} &
			\includegraphics[width=.49\linewidth]{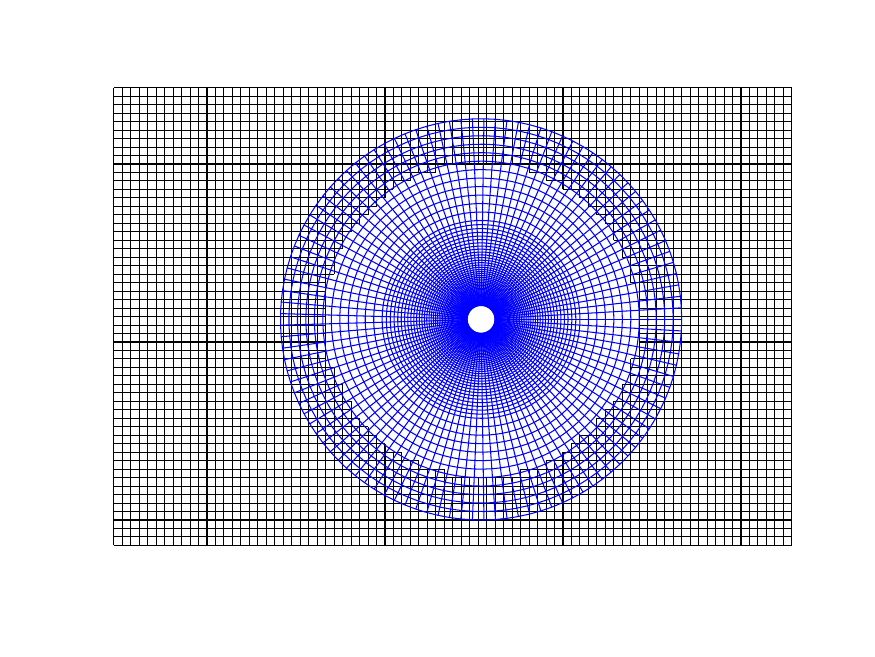}
		\end{tabular}
	\end{center}	
	\caption{Overset configuration for steady and moving cylinder tests. Left: configuration of steady cylinders or the final configuration of moving cylinders. Right: initial configuration of moving cylinders.}
	\label{fig.cylinders_conf}
\end{figure}

When the cylinder does not move (i.e., $\uu_B = \mathbf{0}$), the cylinder has its center $\xx_c$ at the origin of the axes. Consequently, in order to limit possible boundary effects, the inlet and outlet conditions are distant $8D$ and $16D$ in $x-$direction, respectively. The horizontal walls are both at a distance of $8D$ from $\xx_c$ along the $y-$direction. This configuration is shown in Figure \ref{fig.cylinders_conf} (left). On the other hand, when the cylinder moves with a constant horizontal velocity $\uu_B = [-u_B, 0]^\top$, with $u_B > 0$, the just described configuration is always the final configuration at time $t = t_f$. Thus, the initial configuration is obliged to have the center of the cylinder at position $\xx_c = (t_f u_B, 0)$, as sketched in Figure \ref{fig.cylinders_conf} (right). 

The Chimera mesh is defined by a Cartesian background of cell size $\bar{h} = 0.3$. The foreground mesh is built around the cylinder. Its cells have size varying from $0.079$ (close to $\partial \mathcal{C}$) up to $0.3$ (for fringe cells) for a total number of $N_c = 12358$ active cells. For all cases, the fluid velocity is obliged to be equal to $\uu_\infty$ at the inlet and to $\uu_B$ on $\partial \mathcal{C}$. No reflecting conditions are imposed at the outlet \cite{jin1993nonreflecting}, i.e. $[\nabla \uu] \nn = \mathbf{0} $. Free streamline conditions (namely $v = 0$ and $\partial_y u = 0$) close the problem on the horizontal walls. The pressure is strongly put to zero on the outlet. On the remaining boundaries, homogeneous Neumann conditions are imposed, i.e., $\nabla p \cdot \nn = 0$. It is clear that zero inlet velocity and horizontal body velocity $\uu_B = \uu_c$ defines an equivalent fluid dynamics system of a steady cylinder with a fluid entering in the channel with a prescribed nonzero velocity $\uu_\infty = - \uu_c$ at the inlet. For this reason, when the cylinder does not move, the initial fluid velocity is constantly equal to $\uu_\infty = [1, 0]^\top$; when the cylinder moves of velocity $\uu_B = [-1,0]^\top$, the fluid is originally at rest. In both cases, the pressure is zero throughout the entire domain.

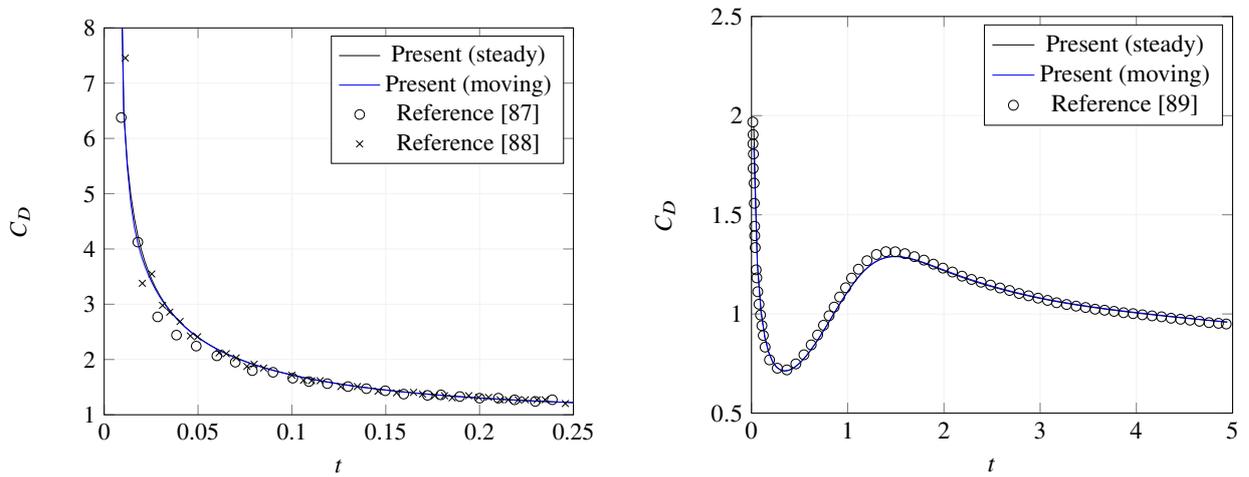
\begin{figure}[!htbp]
	\begin{center}
		\begin{tabular}{cc}
			\resizebox{.49\linewidth}{!}{ \begin{tikzpicture}
				\begin{axis}
					[
					xlabel = {$t$},
					ylabel = {$C_D$},
					xmin = 0, xmax = 0.25,
					ymin = 1, ymax = 8,
					xtick= {0, 0.05, ..., 0.26},
					ytick= {1, 2, ..., 8},
					xticklabels={0,0.05,0.1,0.15,0.2,0.25},
					line/.style={thick},
					grid=both,
					grid style={line width=.1pt, draw=gray!10}]
					
					\addplot[black, solid] table {static_cylinder_Re200.dat}; \addlegendentry{Present (steady)}
					\addplot[blue, solid] table {moving_cylinder_Re200.dat}; \addlegendentry{Present (moving)}
					\addplot[black, mark=o, only marks] table {Koum_Re200.dat}; \addlegendentry{Reference \cite{koumoutsakos1995high} }
					\addplot[black, mark=x, only marks] table {Berg_Re200.dat}; \addlegendentry{Reference \cite{bergmann2014accurate} }
				\end{axis}
			\end{tikzpicture} } &
			\resizebox{.49\linewidth}{!}{ \begin{tikzpicture}
				\begin{axis}
					[
					xlabel = {$t$},
					ylabel = {$C_D$},
					xmin = 0, xmax = 5,
					ymin = 0.5, ymax = 2.5,
					line/.style={thick},
					grid=both,
					grid style={line width=.1pt, draw=gray!10}]
					
					\addplot[black, solid] table {static_cylinder_Re550.dat}; \addlegendentry{Present (steady)}
					\addplot[blue, solid] table {moving_cylinder_Re550.dat}; \addlegendentry{Present (moving)}
					\addplot[black, mark=o, only marks] table {PW2000.dat}; \addlegendentry{Reference \cite{ploumhans2000vortex} }
				\end{axis}
			\end{tikzpicture} }
		\end{tabular}
	\end{center}
	\caption{Time evolution of the drag coefficient $C_D$ for steady and moving cylinders against reference solutions. The Reynolds number is $\Rey = 200$ and $\Rey = 550$ on the left and right, respectively.}
	\label{fig.cylinders_CD}
\end{figure}

We consider different cases with two Reynolds number, that is $\Rey \in \{ 200, 550 \}$. For both viscosity values, we compare the drag coefficients with data from literature. In particular, in Figure \ref{fig.cylinders_CD} (left), we compare the drag coefficient for steady and unsteady cylinder at Reynolds $\Rey = 200$ up to the final time $t_f = 0.25$ with data from \cite{koumoutsakos1995high, bergmann2014accurate}. The same comparison is performed for $\Rey = 550$ with data from \cite{ploumhans2000vortex} in Figure \ref{fig.cylinders_CD} (right) up to the final time $t_f = 5$. In both cases, an excellent matching can be appreciated with the numerical results of our method. 

\subsection{Strouhal number for laminar flow past a cylindrical obstacle}
Let us consider the previously introduced test case of a steady cylinder with Reynolds number $\Rey = 200$. We simulate the laminar flow on a longer time interval of range $t \in [0, 90]$ in order to test the numerical scheme up to the asymptotic regime of this physical phenomenon. We perform a comparison of the average drag coefficient $\bar{C}_D$ and Strouhal number $\St$ at the asymptotic regime. In particular, the Strouhal number is defined as $\St = f_v D / | \uu_c |$, with $f_v$ the frequency of oscillation of the lift coefficient $C_L$ at the asymptotic regime.

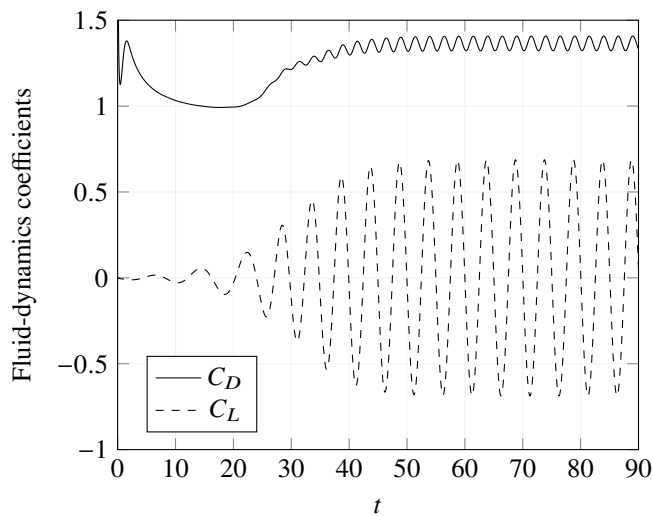
\begin{figure}[!htbp]
	\begin{center}
		\begin{tikzpicture}
			\begin{axis}
				[
				xlabel = {$t$},
				ylabel = {Fluid-dynamics coefficients},
				xmin = 0, xmax = 90,
				ymin = -1, ymax = 1.5,
				xtick= {0, 10, ..., 90},
				line/.style={thick},
				legend style = {at={(axis cs:5,-.9)}, anchor=south west},
				grid=both,
				grid style={line width=.1pt, draw=gray!10}]
				
				\addplot[black, solid] table {CD_Re200long.dat}; \addlegendentry{$C_D$}
				\addplot[black, dashed] table {CL_Re200long.dat}; \addlegendentry{$C_L$}
			\end{axis}
		\end{tikzpicture}
	\end{center}	
	\caption{Drag $C_D$ and lift $C_L$ coefficient as functions of time for the laminar flow over a static cylinder at Reynolds $\Rey = 200$.}
	\label{fig.CDCL_cyl_Re200_long}
\end{figure}

\begin{figure}[!htbp]
	\begin{center}
		\begin{tabular}{c}
			\includegraphics[trim={0 3cm 0 3cm},clip,width=.9\linewidth]{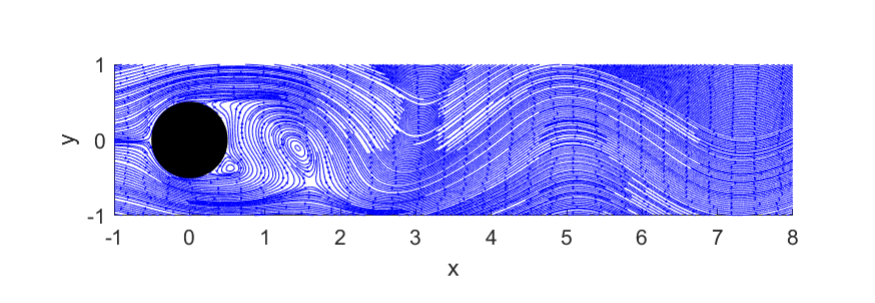}
		\end{tabular}
	\end{center}	
	\caption{Streamlines at the final time $t = 90$ for the laminar flow over a static cylinder at Reynolds $\Rey = 200$.}
	\label{fig.magn_sl_cyl_Re200long}
\end{figure}

The results are collected in Table \ref{tab.steady_Re200} showing a very good matching of our results against the ones from the literature for both the average drag coefficient $\bar{C}_D$ and the Strouhal number $\St$. In Figure \ref{fig.CDCL_cyl_Re200_long} there are the plots of the time evolution of the drag and lift coefficients. In Figure \ref{fig.magn_sl_cyl_Re200long} we report the streamlines at the final time $t = 90$ in a relevant subdomain downstream of the cylinder, which highlights the presence of vortical patterns past the obstacle.

\begin{table}[!htbp]
	\caption{Comparison for the average drag coefficient $\bar{C}_D$ and the Strouhal number $\St$ for the steady cylinder at $\Rey = 200$.}
	\label{tab.steady_Re200}
	\centering
	\begin{tabular}{lll} \toprule
		& $\bar{C}_D$ & $\St$ \\ \midrule
		Present & 1.3661 & 0.1989 \\
		Reference \cite{bergmann2022second} & 1.3430 & 0.1979 \\
		Reference \cite{bergmann2004optimisation} & 1.3900 & 0.1999 \\
		Reference \cite{bergmann2011modeling} & 1.3500 & 0.1980 \\
		Reference \cite{bergmann2014accurate} & 1.4000 & $-$ \\
		Reference \cite{braza1986numerical} & 1.4000 & 0.2000 \\
		Reference \cite{he2000active} & 1.3560 & 0.1978 \\
		Reference \cite{henderson1995details} & 1.3412 & 0.1971 \\ \bottomrule
	\end{tabular}
\end{table}

\subsection{Revolution of a cylinder}
We close the suite of numerical tests by presenting a qualitative analysis of a fluid subject to a rotating cylinder. The computational domain is the channel $[-20, 20] \times [-10, 10]$ with an embedded cylinder of radius $r_c = 0.5$. The center of the cylinder is originally located at $\xx_c(0) = (0,1)$. A counterclockwise rotation around the origin of the axis is imposed with velocity $\dot{\xx}_c(t) = [y_c(t), -x_c(t)]^\top$. At the initial time, the fluid is at rest. The boundary conditions follow the ones of the already introduced test cases with cylinders in Section \ref{sec.cyl}. The background mesh has a cell size $\bar{h} = 0.5$, while the foreground grid is a polar mesh built around the cylinder structure. It is finer near the cylinder (with a cell size $\bar{h} = 0.015$) with fringe cells of size $\bar{h} = 0.4875$. The total number of time-averaged active cells is  $\bar{N}_c = 10620$. The foreground mesh moves according to the displacement of the cylinder, hence obtaining a deformation velocity of $\w = [y, -x]^\top$. For this test case, we consider 10 complete revolutions of the cylinder.

\begin{figure}[!htbp]
	\begin{center}
		\begin{tabular}{ccc}
				\includegraphics[width=.3\linewidth]{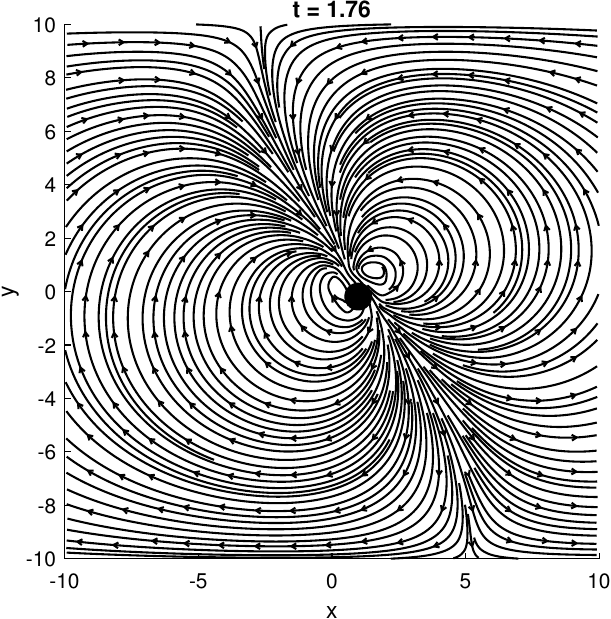} &
				\includegraphics[width=.3\linewidth]{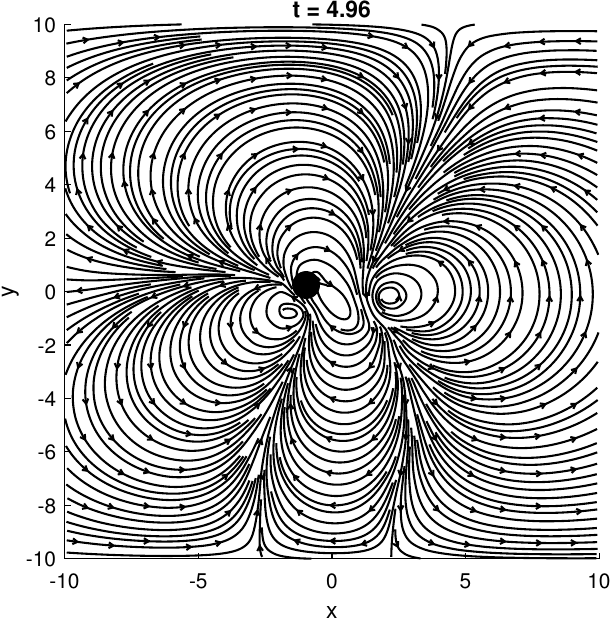} &
				\includegraphics[width=.3\linewidth]{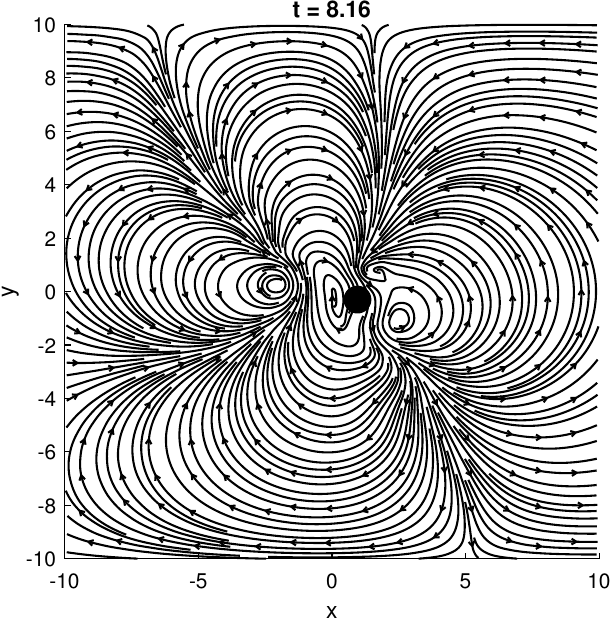} \\
				\includegraphics[width=.3\linewidth]{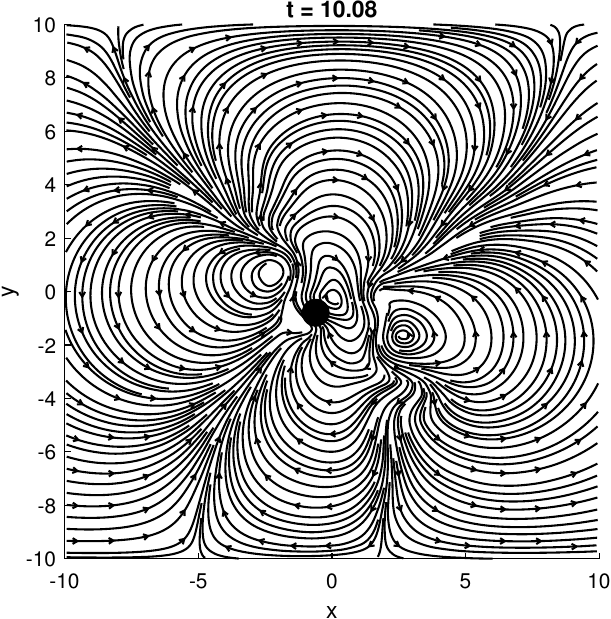} &
				\includegraphics[width=.3\linewidth]{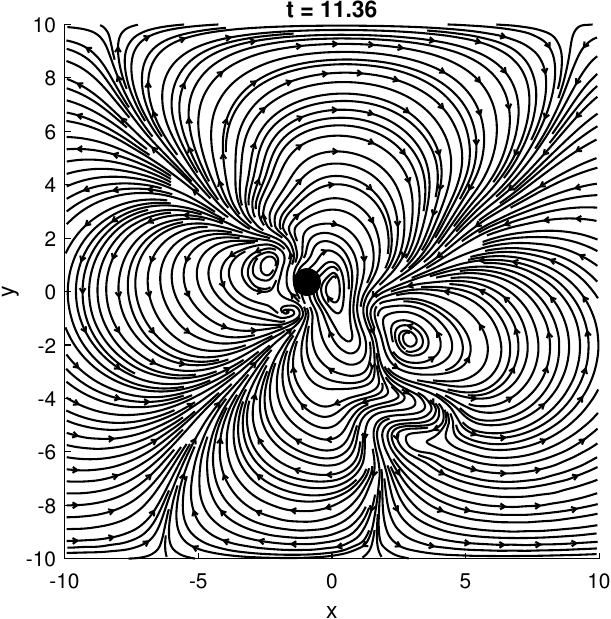} &
				\includegraphics[width=.3\linewidth]{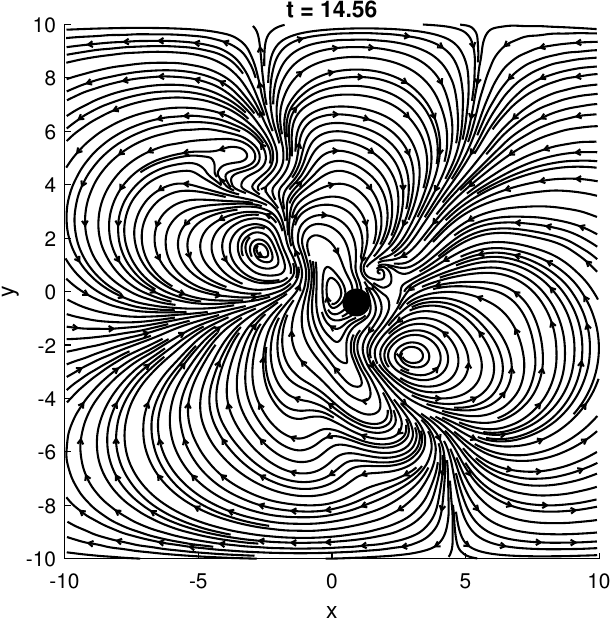} \\
				\includegraphics[width=.3\linewidth]{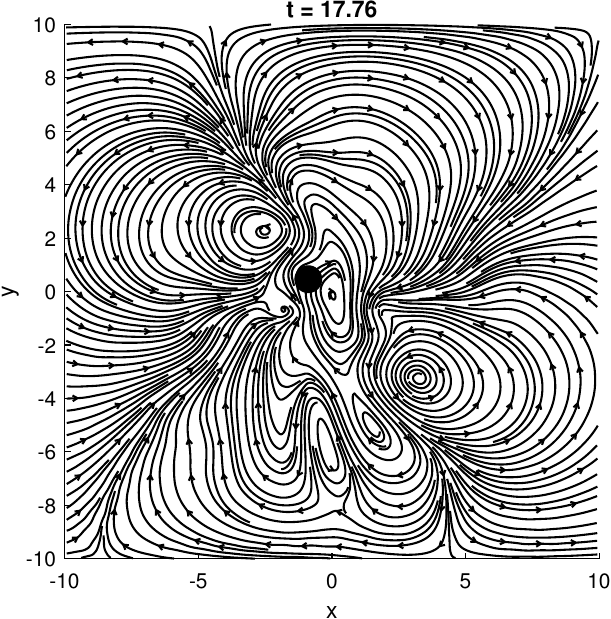} &
				\includegraphics[width=.3\linewidth]{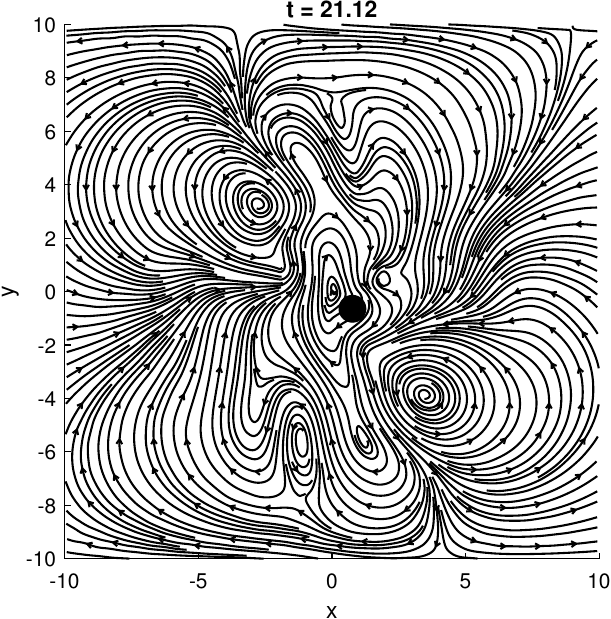} &
				\includegraphics[width=.3\linewidth]{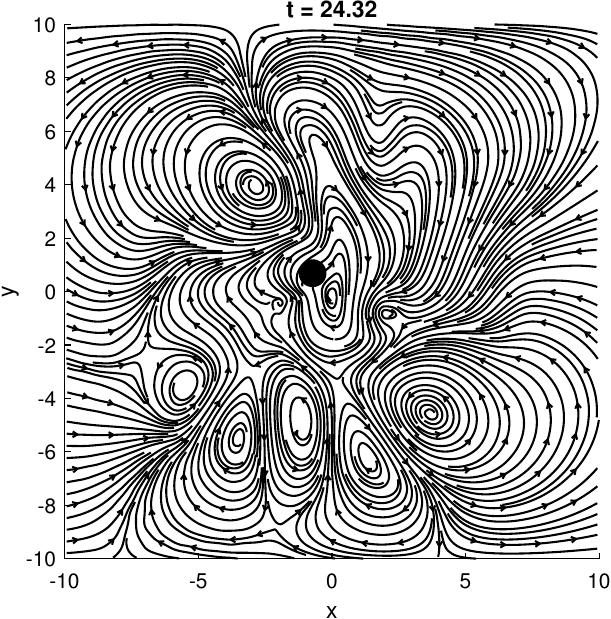} \\
				\includegraphics[width=.3\linewidth]{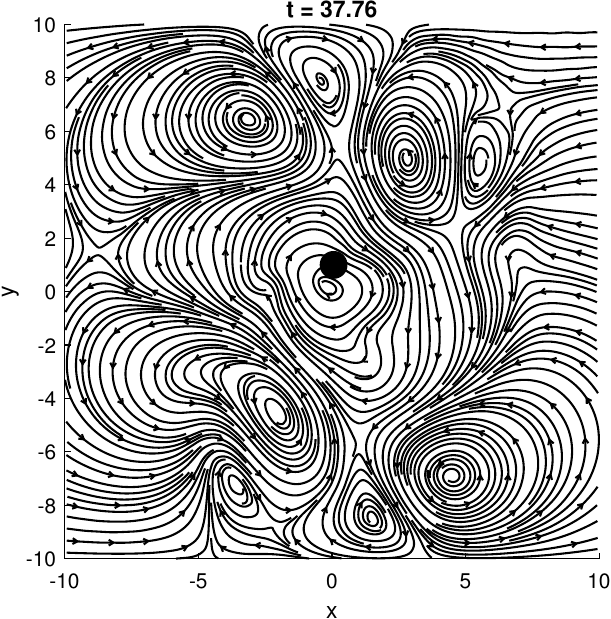} &
				\includegraphics[width=.3\linewidth]{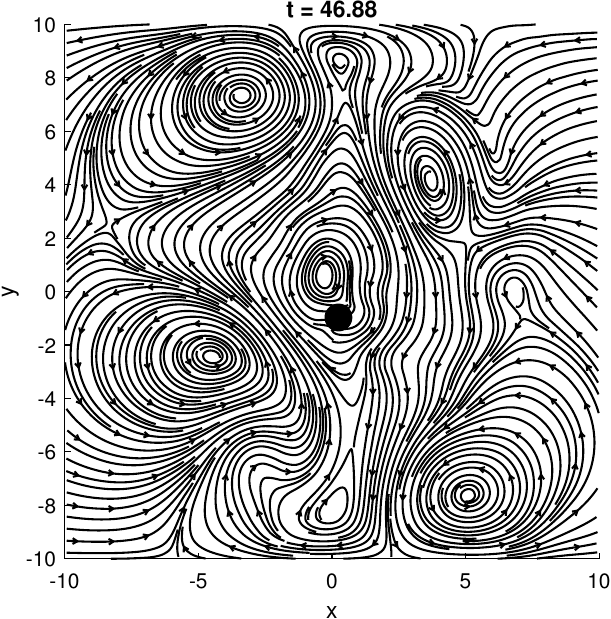} &
				\includegraphics[width=.3\linewidth]{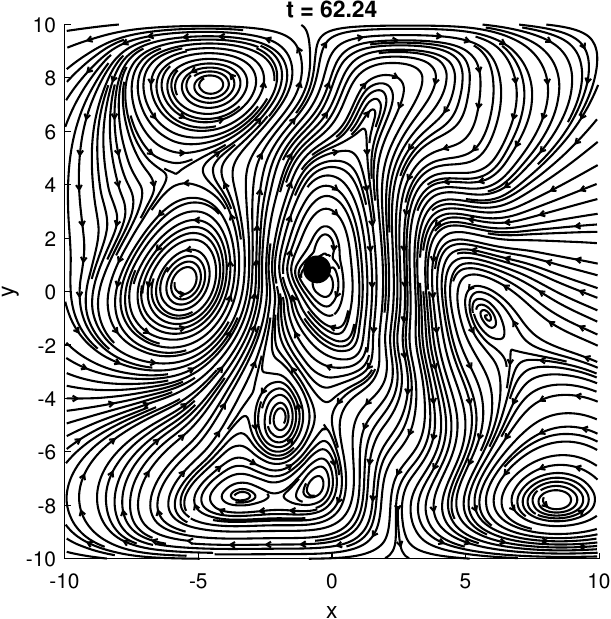} 
		\end{tabular}
	\end{center}
	\caption{Streamlines for the rotating cylinder for different time instances.}
	\label{fig.rotating_cylinder}
\end{figure}

Figure \ref{fig.rotating_cylinder} shows the streamlines of the numerical simulation for 12 significant time instances. From the first times, it is possible to remark vortexes coming off the wake produced by the circular motion of the cylinder. In the meanwhile, a central vortex confined into the trajectory of the center of the cylinder arises. At first, the vortices that detach from the wake are subjected to centrifugal motion, which then pulls them away from the cylinder. Subsequently, their centrifugal motion stops and they begin to circulate counterclockwise in harmony with the motion of the cylinder.  After time $t \simeq 21$, no significant detaching vortexes can be appreciated and the old created vortexes keep rotating with the cylinder (as it is possible to see in the last row of Figure \ref{fig.rotating_cylinder}).   

\section{Conclusions} \label{sec.concl}
The objective of this work has been the design of a second order finite volume numerical method for the incompressible Navier-Stokes equations on moving Chimera meshes. We numerically integrate in time the arising non-autonomous system of partial differential equations through a semi-implicit IMEX scheme. This allows to separate the fast and slow scales of the phenomenon under consideration through a flux splitting technique. The evolution of the frame is encoded in the numerical flux, hence performing integration on active cells of the overset configuration, i.e. on any evolving control volume discretizing the computational domain at the current time level. An incremental fractional-step approach is employed in a context of projection-correction methods in order to ensure the velocity field to be divergence-free. Fringe cells are treated in the donor/receptor paradigm. Therefore, a compact continuous transmission is provided in the numerical approximation of any differential problem. This exploits the continuity of the solution and it is given by an extrapolation of available data in the minimal stencil of the other partition defining a minimal neighborhood with respect to the cell center of fringe cells. Moreover, this allows to transfer information from one block to another in the same time the algebraic system is solved, without creating \textit{ad hoc} discretizations for fringe cells and without exploiting iterative processes (e.g. Schwartz method) in order to ensure the continuity of the solution on the overlapping zone.

The numerical validation demonstrates the second order convergence behavior regardless the evolution of the foreground mesh, even for time-dependent deformations. Furthermore, the new schemes are asymptotically accurate, meaning that the formal order of accuracy is attained independently of the stiffness of the problem under consideration, e.g. independently of the Reynolds number. Moreover, the method is precise at zero-machine if the solution is a polynomial of degree less than or equal to two on moving overset grids. The novel numerical method is also numerically proven to be compliant with the Geometric Conservation Law by performing a free-stream preservation test with different velocities of the Chimera mesh. The presence of the foreground mesh and the overlapping zone does not affect in any sense the solution, as tested for lid-driven cavity flows. Finally, different benchmarks with problems of channels filled of fluids impacting over a cylinder obstacle are presented. They show the accuracy and precision of the novel numerical technique when compared against data from the literature for different Reynolds numbers.

A 3D extension of the method is devised to account for more complex flows. In the future, we plan to exploit this method also for compressible flows where conservation properties are crucial. To approach more realistic scenarios, the adoption of unstructured grids is likely to be pursued. From the viewpoint of the numerical method, the usage of hybrid finite volume/finite element methods for the slow and fast scales is foreseen, in the optic of \cite{boscheri2023new}.  

\section*{Acknowledgments}

WB and MGC received financial support by Fondazione Cariplo and Fondazione CDP (Italy) under the project No. 2022-1895. WB also acknowledges funding from the Italian Ministry of University and Research (MUR) with the PRIN Project 2022 No. 2022N9BM3N. This work was partially carried out at the Institute des Mathématiques de Bordeaux (IMB, Bordeaux-France) during the visiting program of WB and MGC.

\bibliographystyle{elsarticle-num}
\bibliography{biblio}

\end{document}